\documentclass[a4paper,11pt]{amsart}

\usepackage{amsmath}
\usepackage{amscd}
\usepackage{amsxtra}
\usepackage{amssymb}
\usepackage{amsthm}

\newtheorem{thm}{Theorem}[section]
\newtheorem{cor}[thm]{Corollary}

\newtheorem{prop}[thm]{Proposition}
\newtheorem{conj}[thm]{Conjecture}

\newtheorem*{thm*}{Theorem}
\newtheorem*{prop*}{Proposition}
\newtheorem*{cor*}{Corollary}
\newtheorem*{conj*}{Conjecture}

\theoremstyle{definition}
\newtheorem{defn}[thm]{Definition}

\theoremstyle{remark}
\newtheorem{ex}[thm]{Example}
\newtheorem{rem}[thm]{Remark}

\newcommand{\ka}{{\mathcal A}}

\newcommand{\kd}{{\mathcal D}}

\newcommand{\kf}{{\mathcal F}}

\newcommand{\kh}{{\mathcal H}}

\newcommand{\ko}{{\mathcal O}}

\newcommand{\IA}{{\mathbb A}}

\newcommand{\IC}{{\mathbb C}}

\newcommand{\IH}{{\mathbb H}}

\newcommand{\IP}{{\mathbb P}}
\newcommand{\IQ}{{\mathbb Q}}
\newcommand{\IR}{{\mathbb R}}

\newcommand{\IZ}{{\mathbb Z}}

\newcommand{\curly}[1]{\mathcal{#1}}

\newcommand{\DR}{\mathrm{\mathbb R}} 

\newcommand{\eps}{\varepsilon}

\newcommand{\id}{{\rm id}}
\newcommand{\tensor}{\otimes}

\newcommand{\ra}{\rightarrow}
\newcommand{\lra}{\longrightarrow}

\newcommand{\isom}{\cong}

\newcommand{\del}{\partial}
\newcommand{\union}{\cup}

\newcommand{\xlra}[1]{\overset{#1}{\lra}}

\newcommand{\qtext}[1]{\quad\text{#1}\quad}

\newcommand{\Pic}{\mathrm{Pic}}

\newcommand{\Set}[2]{\left\{\, #1 \;|\; #2 \,\right\}}

\newcommand{\HH}{\mathrm{H}}

\newcommand{\<}{\langle}
\renewcommand{\>}{\rangle}

\newcommand{\vsum}{\oplus}

\renewcommand{\HH}{H}

\newcommand{\dual}{\hspace{0cm}^\vee}

\newcommand{\HT}{\tilde{\HH}}

\newcommand{\Hodge}[1]{\mathcal{#1}}

\newcommand{\HOZ}{\Hodge{H}_\IZ}

\newcommand{\HO}{\Hodge{H}}

\newcommand{\QF}{(\,.\,)}
\newcommand{\Q}[2]{(#1.#2)}

\newcommand{\PD}{\mathscr{D}}
\newcommand{\PER}{\mathscr{P}}
\newcommand{\EL}{\tilde{\Lambda}}
\newcommand{\KL}{\Lambda}

\newcommand{\HA}{\HH_A}
\newcommand{\HB}{\HH_B}

\renewcommand{\Re}{Re}
\renewcommand{\Im}{Im}

\usepackage[all,cmtip]{xy}
\usepackage{graphicx}
\usepackage{mathrsfs}
\usepackage{verbatim}

\title{Period- and Mirror-maps for the Quartic K3}
\author{Heinrich Hartmann}
\date{\today}

\begin{document}
\begin{abstract}
  We study in detail mirror symmetry for the quartic K3 surface in
  $\IP^3$ and the mirror family obtained by the orbifold construction.
  As explained by Aspinwall
  and Morrison \cite{AM1997}, mirror symmetry for K3 surfaces can be
  entirely described in terms of Hodge structures.
  \\
  - We give an explicit computation of the Hodge structures and period
  maps for these families of K3 surfaces.
  \\
  - We identify a mirror map, i.e. an isomorphism between the complex
  and symplectic deformation parameters and explicit isomorphisms
  between the Hodge structures at these points.
  \\
  - We show compatibility of our mirror map with the one defined by
  Morrison \cite{Morrison1992} near the point of maximal unipotent
  monodromy.
  \\
  Our results rely on earlier work by Narumiyah--Shiga \cite{NS2001},
  Dolgachev \cite{Dolgachev1996} and Nagura--Sugiyama \cite{NS1995}.
\end{abstract}

\maketitle
\tableofcontents

\section{Introduction}
Let $(X,I,\omega_X)$ be a Calabi--Yau manifold with complex structure
$I$ and chosen K\"ahler form $\omega_X$.  The philosophy of mirror
symmetry says that certain invariants of the {\em complex} manifold
$(X,I)$ should be encoded by the {\em symplectic} structure $\omega_Y$
of a mirror Calabi--Yau $(Y,J,\omega_Y)$ and vice versa.

Following Aspinwall and Morrison \cite{AM1997} (see also
\cite{HuyTrieste2004} and \cite{HG2005}), mirror symmetry for K3
surfaces can be described in terms of Hodge structures. To a K3
surface $(X,I,\omega_X)$ with chosen K\"ahler form $\omega_X$ we
associate two Hodge structures $\HH_A,\HH_B$ on the lattice
$\HH^*(X,\IZ)$. The essential fact is, that $\HH_A$ only depends on
the symplectic form $\omega_X$ and $\HH_B$ only on the complex
structure $I$.  Now $(X,I)$ is said to be mirror dual to
$(Y,\omega_Y)$ in the Hodge theoretic sense if there exists a Hodge
isometry
\[ \HH_B(X,\IZ) \isom \HH_A(Y,\IZ). \]

This definition can be seen as a refinement of Dolgachev's
\cite{Dolgachev1996} notion of mirror symmetry for families of lattice
polarized K3 surfaces (cf. section \ref{Sec:DolComp}).  There are many
examples of mirror dual families of lattice polarized K3 surfaces, 
e.g. \cite{Belcastro}, \cite{Rohsiepe},
\cite{Dolgachev1996}. On the other hand, the author is not aware of an
explicit example of mirror symmetry in the Hodge theoretic sense in
the literature.

\begin{figure}
  \includegraphics[width=10cm]{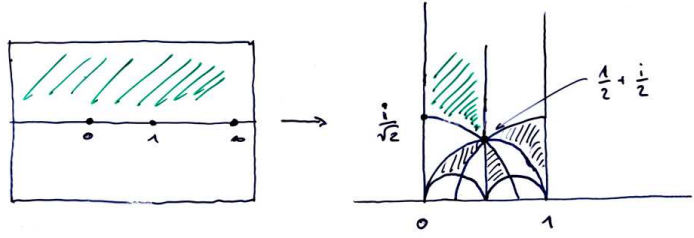}
\caption{Mapping behavior of the mirror map $\psi$ in coordinates $z=1/t^{4}$ and $p$.}
\label{FigureIntro}
\end{figure}

We study the following families of K3 surfaces.
\begin{itemize}
\item Let $Y \subset \IP^3$ be a smooth quartic in $\IP^3$ viewed as a
  symplectic manifold with the symplectic structure given by the
  restriction of the Fubini--Study K\"ahler form $\omega_{FS}$.  
  We introduce a scaling parameter $p \in \IH$ to get a family of
  (complexified) symplectic manifolds $Y_p=(Y,\omega_p), \; \omega_p = p/i \cdot
  \omega_{FS}$ parametrized by the upper half plane.
\item Let $X_t$ be the Dwork family of K3 surfaces, which is
  constructed from the Fermat pencil
  \[ F_t := \{ X_0^4+X_1^4+X_2^4+X_3^4 - 4t X_0 X_1 X_2 X_3 = 0 \}
  \subset \IP^3 \] by taking the quotient with respect to a finite
  group and minimal resolution of singularities.
\end{itemize}
This is the two-dimensional analog to the quintic threefold and its
mirror studied by Candelas et al. \cite{CandelasEtAl}.  
\begin{thm}[Theorem \ref{MMTheorem}, Theorem \ref{mainthm}, Theorem \ref{TriangleThm}]\label{Mainthm}
  The K3 surfaces $X_t$ and $Y_p$ are mirror dual in the Hodge
  theoretic sense if $t$ and $p$ are related by
  \begin{align*}
    exp(2\pi i p) &=
    w+104\,{w}^{2}+15188\,{w}^{3}+2585184\,{w}^{4}+480222434\,{w}^{5}
    +\dots
  \end{align*}
  where $w:=1/(4 t)^{4}$. A closed expression as ratio of
  hypergeometric functions is given in section \ref{NSSols}.

  The multi-valued map $\psi: z \mapsto p(z), z=1/t^4$ determined by
  this equation is a Schwarz triangle function which maps the upper
  half plane to the hyperbolic triangle with vertices
  $(\infty,\frac{i}{\sqrt{2}},\frac{1+i}{2})$ and interior angles
  $(0,\pi/2,\pi/4)$, as pictured in Figure \ref{FigureIntro}.
\end{thm} 

The proof relies heavily on earlier work by Narumiyah and Shiga
\cite{NS2001}, Dolgachev \cite{Dolgachev1996} and Nagua and Sugiyama
\cite{NS1995}.  We proceed in three main steps:  First, we use a
theorem of Narumiyah and Shiga which provides us with the required
cycles and a description of the topological monodromy of the family.
Then we consider the Picard--Fuchs differential equation which is
satisfied by the period integrals.  We derive a criterion for a set of
solutions to be the coefficients of the period map.  In a third step
we construct solutions to this differential equation which match this
criterion. Here we use the work of Nagura and Sugiyama. The relation
to Schwarz triangle function appears also in \cite[Thm. 6.1]{NS2001}.

The function in Theorem \ref{Mainthm} was also considered by Lian and
Yau \cite{LY1996} (see Remark \ref{Comparison}).  There it was noted
that the inverse function $z(p)$ is a modular form with integral
Fourier expansion which is related to the Thompson series for the
Griess--Fischer (``monster'') group.  See also the exposition by
Verrill and Yui in \cite{VY1998}.

Our motivation for studying this specific family stems from a theorem
of Seidel.  Recall that the homological mirror symmetry conjecture due
to Kontsevich \cite{HMS} states that if $(X,I)$ is mirror dual to
$(Y,\omega_Y)$ then then there is an exact equivalence of triangulated
categories
\[ \kd^b(Coh(X,I)) \isom \kd^\pi(Fuk(Y,\omega_Y)). \]

There are only a few cases where such an equivalence is known to hold. 
One example was provided by Seidel. He proves homological mirror symmetry 
for the pair of K3 surfaces considered above.

\begin{thm*}[Seidel \cite{Seidel2003}]
  If the family $X_t$ is viewed as a K3 surface $X$ over the Novikov field $\Lambda_{\IQ}(1/t)$, 
  which is the algebraic closure of the field of formal Laurent series $\IC((1/t))$,
  then there is an isomorphism $\psi: \Lambda_{\IQ}(1/t) \isom \Lambda_{\IQ}(q)$ and an equivalence of 
  triangulated $\Lambda_{\IQ}(q)$-linear categories
  \[  \psi_*\kd^b(Coh(X)) \isom \kd^\pi(Fuk(Y)). \]
\end{thm*}

Unfortunately, the isomorphism $\psi$ has not yet been determined explicitly. 
Geometrically it describes the dependence of the symplectic volume $p$ of the quartic from
the deformation parameter $t$ of the complex structure on $X$.
Thus our mirror map $\psi$ in Theorem \ref{Mainthm} provides a conjectural candidate for this 
isomorphism. 

On the way to proving Theorem \ref{Mainthm} we also give an explicit
calculation of the classical period map for the Dwork family.
Consider a non-zero holomorphic two-form $\Omega \in H^{2,0}(X)$ and a
basis of two-dimensional cycles $\Gamma_i \in \HH_2(X,\IZ)\isom
\IZ^{22}$. By the global Torelli theorem, the complex structure on $X$
is determined by the period integrals $
(\int_{\Gamma_1}\Omega,\dots,\int_{\Gamma_{22}} \Omega) $ and the
intersection numbers $\Gamma_i.\Gamma_j$.

\begin{thm}[Theorem \ref{mainthm},  Remark
  \ref{PeriodThmProof}]\label{PeriodThm}
  For $t \in \IC$ near  $t_0 = i / \sqrt{2}$, 
  there are explicit bases $\Gamma_{i}(t) \in \HH_2(X_t,\IZ),
  i=1,\dots,22$ and holomorphic two forms $\Omega_t$ on the Dwork
  family $X_t$ such that the period integrals are given by
  \[ (\int_{\Gamma_1(t)}\Omega_t,\dots,\int_{\Gamma_{22}(t)} \Omega_t)
  = (4 p(t),2p(t)^2,-1,p(t),0,\dots,0), \] where $p(t)=\psi(1/t^4)$ is
  the function introduced in Theorem \ref{Mainthm}.
\end{thm}

{\bfseries Acknowledgments.}
This work is part of my PhD thesis written under the supervision of Prof. Daniel Huybrechts
to whom I owe much gratitude for his constant support and encouragement.
I thank Duco van Straten for explaining to me much about
hypergeometric functions and Picard--Fuchs equations.


\section{Mirror symmetry for K3 surfaces} \label{Framework} In this
section we summarize Aspinwall and Morrison's description
\cite{AM1997} of mirror symmetry for K3 surfaces in terms of Hodge
structures. Their constructions have been generalized to higher
dimensional hyperk\"ahler manifolds by Huybrechts in
\cite{HuyTrieste2004} and \cite{HG2005}. 

\subsection{The classical Hodge structure of a complex K3 surface}
Recall from \cite{Beauville1985} that a K3 surface is a
two-dimensional connected complex manifold $X$ with trivial canonical
bundle $\Omega^2_X \isom \ko_X$ and $H^1(X,\ko_X)=0$.

The second cohomology $\HH^2(X,\IZ)$ endowed with the cup-product
pairing $\Q{a}{b}=\int a \cup b$ is an even, unimodular lattice of
rank 22 isomorphic to the K3 lattice $\Lambda := 2 E_8(-1) \vsum 3 U$.
The group  $\HH^{2,0}(X)=\HH^0(X,\Omega^2_X)$ is spanned by a the class of a 
holomorphic two form $\Omega$ which is nowhere vanishing. This class satisfies
the properties
\[ \Q{\Omega}{\Omega} = 0\, \qquad \Q{\Omega}{\bar{\Omega}} > 0. \]

\begin{rem}\label{h20_hodge}
  The Hodge structure on $H^2(X,\IZ)$ is completely determined by the subspace $\HH^{2,0}(X)
  \subset \HH^2(X,\IC)$. Indeed, we have 
  \begin{align*}  \HH^{0,2}(X) = \overline{\HH^{2,0}(X)}
    \qtext{and} \HH^{1,1}(X)=(\HH^{2,0}(X)\vsum
    \HH^{0,2}(X))^\perp. 
  \end{align*}
\end{rem}

The global Torelli theorem states that a K3 surface is determined up
to isomorphy, by it's Hodge structure.
\begin{thm}[Piatetski-Shapiro--Shafarevich,
  Burns--Rapoport] \label{GlobalTorelli} Two K3 surfaces $X,X'$ are
  isomorphic if and only if there is a Hodge isometry $ \HH^2(X,\IZ)
  \isom \HH^2(X',\IZ)$.
\end{thm}

\subsection{CFT-Hodge structures of complex K3 surfaces}
There is another weight-two Hodge structure associated to a K3
surface, which plays an important role for mirror symmetry.

Define the {\em Mukai pairing} on the total cohomology $\HH^*(X,\IZ)$ by
\begin{align}\label{MukaiPairing} \Q{(a_0,a_2,a_4)}{(b_0,b_2,b_4)}:=
  \int a_2 \cup b_2 - a_0 \cup b_4 - a_4 \cup b_2.
\end{align}
We denote this lattice by $\HT(X,\IZ)$. It is an even, unimodular
lattice of rank $24$ and signature $(4,20)$ isomorphic to the enlarged
K3 lattice $\tilde{\Lambda}:=\Lambda \oplus U$.

We define a weight-two Hodge structure on $\HT(X,\IZ)$ by setting
$\HB^{2,0}(X)=H^{2,0}(X)$ and using the construction in Remark
\ref{h20_hodge}.  Note that \[ \HB^{1,1}(X)=\HH^0(X,\IC) \vsum
\HH^{1,1}(X) \vsum \HH^4(X,\IC).\] We call
$\HH_B(X,\IZ)=(\HT(X,\IZ),\QF,\HB^{p,q}(X))$ the {\em B-model Hodge
  structure} of $X$.  The name is motivated by the statement in
\cite{AM1997}, that the ``B-model conformal field theory'' associated
to $X$ is uniquely determined by $\HB(X,\IZ)$.

One very important occurrence of this Hodge structure is the following theorem.
\begin{thm}[Derived global Torelli; Orlov \cite{Orlov1997}]
  Two projective K3 surfaces $X,X'$ have equivalent derived categories
  $\kd^b(X) \isom \kd^b(X')$ if and only if there exists a Hodge
  isometry $\HB(X,\IZ) \isom \HB(X',\IZ)$.
\end{thm}

\subsection{CFT-Hodge structures of symplectic K3 surfaces}\label{ACFTHodge}
Every K\"ahler form $\omega$ on a K3 surface $X$ defines a symplectic
structure on the underlying differentiable manifold.  In this section
we will associate a Hodge structure to this symplectic manifold.
Moreover, we shall allow twists by a so called B-field $\beta \in
\HH^2(X,\IR)$ to get a complexified version.

Given $\omega$ and $\beta$ we define the following class of mixed,
even degree
\begin{align}\label{mho} 
  \mho = \exp(i\omega + \beta) = (1,i\omega + \beta,(i\omega +
  \beta)^2/2 ) \in \HT(X,\IC).
\end{align}
This class enjoys formally the same properties as $\Omega \in
\HH^{2,0}(X)$ above:
\[ \Q{\mho}{\mho} = 0, \quad \Q{\mho}{\overline{\mho}}> 0\] with respect to
the Mukai-pairing.  Hence, we can define a Hodge structure
$\HA(X,\IZ)$ on $\HT(X,\IZ)$ by demanding $\HA^{2,0}(X):=\IC\,\mho$
via Remark \ref{h20_hodge}.

We call $\HA(X,\IZ)$ the {\em A-model Hodge structure} of
$(X,\omega,\beta)$. Again, the name is motivated by the statement in
\cite{AM1997}, that the ``A-model conformal field theory''
associated to $(X,\omega,\beta)$ is uniquely determined by
$\HA(X,\IZ)$.

\subsection{Mirror symmetries}
Two Calabi--Yau manifolds $X,Y$ form a mirror pair if the B-model
conformal field theory associated to $X$ is isomorphic to the
A-model conformal field theory associated to $Y$. This motivates
the following definition.

\begin{defn}
  A complex K3 surface $X$ with holomorphic two-form $\Omega$ and a
  symplectic K3 surface $Y$ with complexified K\"ahler form
  $\mho=\exp(i \omega + \beta)$ form a {\em mirror pair} if there
  exists a Hodge isometry
  \[ \HH_B(X,\IZ) \isom \HH_A(Y,\IZ). \]
\end{defn}

Thus a naive translation of Kontsevich's homological mirror conjecture
reads as follows.

\begin{conj}
  Let $X$ be a K3 surface with holomorphic two-form $\Omega$ and 
  $Y$ a K3 surface with K\"ahler form $\omega$.
  Then there is an exact equivalence of triangulated categories
  \[ \kd^b(Coh(X)) \isom \kd^\pi(Fuk(Y)) \]
  if and only if there is a Hodge isometry
  $ \HH_B(X,\IZ) \isom \HH_A(Y,\IZ). $
\end{conj}

Note that this is perfectly consistent with Orlov's derived global Torelli theorem.

\subsection{Relation to mirror symmetry for lattice polarized K3
  surfaces}\label{Sec:DolComp}
In this subsection we compare the Hodge theoretic notion of mirror
symmetry to Dolgachev's version for families of lattice polarized K3 surfaces
\cite{Dolgachev1996}.  See also \cite[Sec. 7.1]{HuyTrieste2004} and
\cite[Sec. 2]{Rohsiepe}.

Let $M \subset \KL$ be a primitive sublattice.  A {\em $M$-polarized
  K3 surface} is a K3 surface $X$ together with a primitive embedding
$i:M \ra Pic(X)$. We call $(X,i)$ pseudo-ample polarized if $i(M)$
contains a numerically effective class of positive self intersection.

Assume that $M$ has the property, that for any two primitive
embeddings $i_1,i_2: M \ra \KL$ there is an isometry $g \in O(\KL)$
such that $i_2=g \circ i_1$.  Then, there is a coarse moduli space
$K_M$ of pseudo-ample $M$-polarized K3 surfaces.

Fix a splitting $M^\perp = U \vsum \hat{M}$. The above condition
ensures, that the isomorphism class of $\hat{M}$ is independent of
this choice.

\begin{defn}
  The \em{mirror moduli space} of $K_M$ is $K_{\hat{M}}$.
\end{defn}

Symplectic structures on a K3 surface $Y$ in $K_{\hat{M}}$ correspond
to points of the mirror moduli space $K_{M}$ in the following way:

Let $(Y,j) \in K_{\hat{M}}$ be an $\hat{M}$-polarized K3 surface with
a marking, i.e. an isometry $n: H^2(Y,\IZ) \ra \KL$, such that $j=n^{-1}|_{\hat{M}}$.  Let
$\omega + i \beta \in H^2(Y,\IC)$ be a complexified symplectic
structure on $Y$, which is compatible with the $\hat{M}$-polarization,
i.e. $\omega + i \beta \in j(\hat{M})_{\IC}$.
Denote by $\mho = exp(i \omega \vsum \beta)$ be the associated period vector.

The chosen splitting $M^{\perp} = \hat{M} \vsum U$ determines an isometry 
of $\xi \in O(\HT(Y,\IZ))$ which interchanges the hyperbolic plane $n(U)$ 
with $H^0(Y,\IZ) \vsum H^4(Y,\IZ)$ and leaves the orthogonal 
complement fixed.

By construction the vector $\Omega := \xi (\mho)$ lies in $n(U)_{\IC} \vsum
\hat{M}_{\IC} \subset H^2(Y,\IC)$.  Note that $\Q{\Omega}{\Omega}=0$ and
$\Q{\Omega}{\bar{\Omega}}>0$.  Hence, by the surjectivity of the
period map \cite[Exp. X]{Beauville1985}, there exists a complex K3 surface $X$ and a isometry $g:H^2(Y,\IZ)
\ra H^2(X,\IZ)$ that maps $\Omega$ into $H^{2,0}(X)$.
Extend $g$ to an isometry of Mukai lattices $\tilde{g}$, then
\[ \tilde{g} \circ \xi: H_A(Y,\IZ) \lra H_B(X,\IZ) \]
is an Hodge isometry.
Moreover, the marking of $Y$ induces an $M$-polarization of $X$ via
\[ i: M \subset \KL \xlra{n} H^2(Y,\IZ) \xlra{g} H^2(X,\IZ). \]
This means $(X,i)$ lies in the mirror moduli space $K_M$.

Conversely, if $\Omega \in H^{2.0}(X)$ is the period vector of a
marked $M$-polarized K3 surface, then $\mho=\xi(\Omega)$ lies
in $H^0(X,\IC) \vsum \hat{M}_{\IC} \vsum H^4(X,\IC)$. Hence $\mho$ is
of the form
\[ \mho = a \, exp(i \omega + \beta)\] for some $\omega,\beta \in M_{\IR}$,
$a \in \IC^*$. 
Indeed, write $\mho=(a,c,b)$ with respect to the above decomposition, then $-2ab+c^2=0$ since $\mho^2=0$.
Therefore $a \neq 0$ and we can set $i \omega + \beta := c/a \in M_{\IC}$.

Note that $\omega^2>0$ since $\mho.\bar{\mho} > 0$.
Now assume, that $\omega$ is represented by a symplectic
form, then $i \omega + \beta$ defines a 
complexified symplectic structure on $Y=X$ such that
  \[ \HB(X,\IZ) \isom \HA(Y,\IZ). \]


\newcommand{\FAM}[1]{#1}

\subsection{Period domains}
In order to compare Hodge structures on different manifolds, it is
convenient to introduce the period domains classifying Hodge structures.

Let $(L,\QF)$ be a lattice. The {\em period domain} associated to $L$
is the complex manifold
\[ \PD(L) := \Set{  [\Omega] \in \IP(L \tensor \IC) }{
  \Q{\Omega}{\Omega} = 0, \; \Q{\Omega}{\bar{\Omega}} > 0 }. \] 
The orthogonal group $O(L,\QF)$ acts on $\PD(L)$ from the left.

The period domain carries a tautological variation of Hodge structures
on the constant local system $L$.  Indeed, the holomorphic vector
bundle $L \tensor \ko_{\PD(L)}$ has a tautological sub-vector bundle
$\kf^2$ with fiber $ \IC\, \Omega \subset L \tensor \IC$
over a point $[\Omega] \in \PD(L)$.
The Hodge filtration is determined by $\kf^2$ via
\begin{align}\label{f2_hodge} 
  \kf^2 \subset \kf^1:=(\kf^2)^\perp \subset L \tensor \ko_{\PD(L)} .
\end{align}

\subsection{Periods of marked complex K3 surfaces}
Let $\pi:\FAM{X} \ra B$ be a smooth family of K3 surfaces.  We have a local system
\[ \HOZ = \DR^2\pi_* \underline{\IZ}_\FAM{X} \] 
on $B$ with stalks isomorphic to the cohomology $H^2(X_t,\IZ)$ of the fiber
$X_t=\pi^{-1}(\{t\})$.  It carries a quadratic form $\QF: \HOZ \tensor
\HOZ \ra \HOZ$ and a holomorphic filtration 
 \[ \kf^2=\pi_* \Omega^2_{X/B} \subset \kf^1:= (\kf^2)^{\perp} \subset \HO:=\HOZ \tensor \ko_B   \]
restricting fiber wise to the cup product pairing and the Hodge
filtration on $\HH^2(X_t,\IC)$,
respectively.

Suppose now, that the local system $\HOZ$ is trivial, and we have
chosen a {\em marking}, i.e. an isometric trivialization
$ m: \HOZ \ra \KL \tensor \underline{\IZ}_B. $
We can transfer the Hodge filtration on $\HOZ$ to the constant 
system $\KL$ via $m$ and get a unique map to the period domain
\[ \PER(\kf^*,m): B \lra \PD(\KL) \] 
with the property that the pull-back of the tautological variation of
Hodge structures agrees with $m(\kf^*)$ as Hodge structures on $\KL
\tensor \underline{\IZ}_B$.
If $\Omega$ is a local section of $\kf^2$, then the period map is explicitly given by
\[ \PER(\kf^*,m)(t) = [m(\Omega(t))] \in \PD(\KL) \] for $t \in B$. 

\subsection{CFT-Periods of marked complex K3 surfaces}
In the same way, we define the periods of the enlarged Hodge structures.
We endow the local system
\[ \tilde{\HO}_\IZ := \DR^*\pi_* \underline{\IZ}_\FAM{X} = 
\DR^0\pi_* \underline{\IZ}_\FAM{X} \vsum \DR^2\pi_* \underline{\IZ}_\FAM{X} \vsum \DR^4\pi_* \underline{\IZ}_\FAM{X} \]
with the Mukai pairing defined by the same formula (\ref{MukaiPairing}) as above.
The associated holomorphic vector bundle $\tilde{\HO}= \tilde{\HO}_\IZ \tensor \ko_B$ carries the 
{\em B-model Hodge filtration}
\[ \kf_B^2 := \pi_* \Omega^2_{\FAM{X}/B} \subset
\kf_B^1:=(\kf_B^2)^\perp \subset \tilde{\HO}. \]

For every marking $\tilde{m}: \tilde{\HO}_\IZ \ra \EL \tensor \underline{\IZ}_B$
of this enlarged local system, we get an associated {\em B-model period map}
\[ \PER_B(\kf^*_B,\tilde{m}): B \lra \PD(\EL). \]

\begin{rem}\label{InducedMarking}
A marking $m$ of $\HOZ$ determines a marking of $\tilde{\HO}_\IZ$ by the following convention.
There are canonical trivializing sections $1 \in \DR^0\pi_* \underline{\IZ}_\FAM{X}$ 
and $or \in \DR^4\pi_* \underline{\IZ}_\FAM{X}$, satisfying $\Q{1}{or}=-1$ with respect to the Mukai pairing.
Let $e,f$ be the standard basis of $U$ with intersections $\Q{e}{f}=1,\Q{e}{e}=\Q{f}{f}=0$.
Then the map 
 \[ m_0 :\DR^0\pi_* \underline{\IZ}_\FAM{X} \vsum \DR^4\pi_* \underline{\IZ}_\FAM{X} \lra U \tensor  \underline{\IZ}_B, 
 \quad 1 \mapsto e, \; or \mapsto -f \]
is an orthogonal isomorphism of local systems and the map
\[ \tilde{m}: = m \vsum m_0: \tilde{\HO}_\IZ = \HOZ \vsum (\DR^0\pi_* \underline{\IZ}_\FAM{X} \vsum \DR^4\pi_* \underline{\IZ}_\FAM{X}) \lra 
(\KL \vsum U)\tensor \underline{\IZ}_B = \EL \tensor \underline{\IZ}_B. \]
defines a marking of $\tilde{\HO}_\IZ$.

\end{rem}

\subsection{CFT-Periods of marked symplectic K3 surfaces}
Let $\pi:\FAM{X} \ra B$ be a family of K3 surfaces, and $\omega \in
\HH^0(B,\pi_*\ka^2_{\FAM{X}/B})$ a $d_{\FAM{X}/B}$-closed
two-form, that restricts to a K\"ahler form on each fiber $X_t$.  The
form $\omega$ determines a global section of
\[ \HO_\infty = (\DR^2 \pi_* \underline{\IZ}_\FAM{X}) \tensor \mathscr{C}_B^\infty(\IC)  
= \DR^2 \pi_* (\ka^*_{\FAM{X}/B}) = \kh^2(\pi_* \ka^*_{\FAM{X}/B}).\]
Analogously, a closed form $\beta \in
\HH^0(B,\pi_*\ka^2_{\FAM{X}/B})$ gives a section $\beta \in
\HH^0(B,\HO_\infty)$.  Given $\omega$ and $\beta$ we define a section
\[ \mho=\exp(i \omega + \beta) \in \HH^0(B,\tilde{\HO}_\infty), 
\quad \tilde{\HO}_\infty= \DR^* \pi_* \IZ \tensor
\mathscr{C}_B^\infty(\IC) \]
by the same formula (\ref{mho}) used in the point-wise definition of
$\mho$.  We set the {\em A-model Hodge filtration} to be the sequence
of $\mathscr{C}^\infty$-vector bundles 
\[ \kf^2_A :=  \mathscr{C}_B^\infty(\IC)\, \mho \subset
\kf^1_A:=(\kf^2_A)^\perp \subset \tilde{\HO}_\infty.
\]

In the same way as above, every marking
$ \tilde{m}: \tilde{\HO}_\IZ \ra \EL \tensor \underline{\IZ}_B$
determines an {\em A-model period map}
\[ \PER_A(\kf_A^*,\tilde{m}): B \lra \PD(\EL) \]
which is a morphism of $\mathscr{C}^\infty$-manifolds.

\begin{ex} \label{HolSymplPeriod} Given $d_\FAM{X}$-closed two-forms
  $\omega, \beta \in \ka^2_{\FAM{X}}$ on $\FAM{X}$, we get
  $d_{X/B}$-closed relative two-forms via the canonical projection
  $\ka^*_\FAM{X} \ra \ka^*_{\FAM{X}/B}$.
  In this case, the map $B \ni t \mapsto \exp(i\omega(t) + \beta(t))
  \in \HT(X_t,\IC)$ factors through the pull-back
  \[ i^*: \HH^*(\FAM{X},\IC) \lra \HH^*(X_t,\IC) \]
  along the inclusion $i: X_t \ra \FAM{X}$. As this map is already defined on $H^2(\_,\IZ)$ the 
  associated period map is constant.

  We can extend this example a bit further. Let $\omega$ be a constant K\"ahler form as above
  and $f: B \ra \IH$ a holomorphic function to the upper half-plane.
  The form $f \omega  = i \Im(f) \omega + \Re(f) \omega$ is $d_{\FAM{X}/B}$-closed 
  and satisfies $\Q{\exp(f \omega)}{\overline{\exp(f \omega)}}>0$. Hence we get a period map
  \[ \PER_A(\kf_A^*,\tilde{m})(t) = [\tilde{m}( \exp( f(t) \omega ))] \in \PD(\EL) \]
  for $t \in B$, which is easily seen to be holomorphic.
\end{ex}

\subsection{Mirror symmetry for families}\label{MSFamily}
Let $\pi: \FAM{X} \ra B$ be a family of complex K3 surfaces with marking 
and $\rho: \FAM{Y} \ra C$ a family of K3 surfaces with marking and
chosen relative complexified K\"ahler form $i \omega + \beta$.

\begin{defn} 
  A {\em mirror symmetry} between $\FAM{X}$ and $\FAM{Y}$ consists
  of an orthogonal transformation $g \in O(\EL)$ called {\em global
    mirror map} and an \'etale, surjective morphism $\psi: C \ra B$
  called {\em geometric mirror map}\footnote{ We think of $\psi$ as a
    multi-valued isomorphism: In practice the period map for
    $\FAM{X}$ is only well defined after base-change to a covering
    space $\tilde{B} \ra B$. Moreover, $\psi$ induces an isomorphism
    between the universal covering spaces of $C$ and $B$.  }  such that the following
  diagram is commutative.
  \begin{center}
    \begin{minipage}{5cm} 
      \xymatrix{ 
        C \ar[rr]^{\PER_A(\FAM{Y})}\ar[d]_{\psi} & &  \PD(\EL)\ar[d]^g \\
        B \ar[rr]^{\PER_B(\FAM{X})} & & \PD(\EL) \\
      }
    \end{minipage}
  \end{center}
  In particular for every point $s \in C$ we have a mirror pair
  \[ \HH_A(Y_s,\IZ) \isom \HH_B(X_{\psi(s)},\IZ). \]
\end{defn}

\begin{rem}
  A typical global mirror map will exchange the hyperbolic plane
  $\HH^0(X_s,\IZ) \vsum \HH^4(X_s,\IZ)$ with  a hyperbolic
  plane inside $\HH^2(X_s,\IZ)$ as in subsection \ref{Sec:DolComp}. We
  will see, that this happens in our case, too. Examples for
  other mirror maps can be found in 
  \cite[Sec. 6.4]{HuyTrieste2004}.

  Note that, if the markings of $\FAM{X}$ and $\FAM{Y}$ are both induced by a
  marking of the second cohomology local system as in Remark
  \ref{InducedMarking} then $g$ can never be the identity.
  Indeed, we always have
  \[ \HH^{2,0}_B(X_s) \perp (\HH^0(X_s,\IZ) \vsum \HH^4(X_s,\IZ)) \]
  but never $\HH^{2,0}_A(Y_t) \perp (\HH^0(Y_t,\IZ) \vsum
  \HH^4(Y_t,\IZ))$ since $ \Q{\exp(i \omega + \beta)}{or} = -1.$
\end{rem}


\section{Period map for the quartic}\label{QuarticSection}
Since the calculation of the period map for the symplectic quartic is
much easier than for the Dwork family, we begin with this
construction.

A smooth quartic in $Y \subset \IP^3$ inherits a symplectic structure
from $\IP^3$ by restricting the Fubini--Study K\"ahler form
$\omega_{FS}$. A classical result of Moser \cite{Moser} shows that all
quartics are symplectomorphic.

\begin{prop}
  For all primitive $h \in \KL$ with $\<h,h\> = 4$ there exists a 
  marking $m: \HH^2(Y,\IZ) \ra \KL$ such that $m([\omega])=h$.
\end{prop}

\begin{proof} 
  Recall that $[\omega] = [\omega_{FS}|_Y] \in \HH^2(Y,\IR)$ is an
  integral class and satisfies $\int_Y \omega^2=4$.  Moreover
  $[\omega]$ is primitive since there is an integral class $l$,
  represented by a line on $Y$, with $l.h = 1$.  Let $n: \HH^2(Y,\IZ)
  \ra \KL$ be an arbitrary marking.  We can apply a theorem of
  Nikulin, which we state in full generality below (\ref{NIK_THM}), to
  get an isometry of $\HH^2(Y,\IZ)$ that maps $[\omega]$ to the
  primitive vector $n(h)$ of square $4$.
\end{proof}

Fix a quartic $Y$ with symplectic form $\omega$.  
Scaling the symplectic form by $\lambda \in \IR_{>0}$ and
introducing a B-field $\beta = \mu\, \omega \in \HH^2(X,\IR), \mu \in
\IR$.  We get a family of
complexified symplectic manifolds
\[ \rho:\curly{Y} \lra \IH \] 
with fiber $(Y,\mho=\exp(i p \omega))$ over a point $p=i \lambda + \mu \in \IH$. 

Since the family is topologically trivial, the marking $m$ of $Y$
constructed above extends to a marking of $\rho:\curly{Y} \ra \IH$, which
induces an enlarged marking
\[ \tilde{m}: R^*\rho_* \underline{\IZ}_\curly{Y} \lra \underline{\IZ}_\IH \tensor \EL. \]
by the procedure explained in Remark \ref{InducedMarking}.

\begin{prop} \label{APeriodDomain}
  The A-model period map of the family $\rho:\curly{Y} \ra \IH$
  \[ \PER_A(\kf^*_A,\tilde{m}) : \IH \lra \PD(\EL)  \]
  is holomorphic and induces an isomorphism of $\IH$ onto a
  connected component $\PD(\< h \> \vsum U)^+$ of
  \[ \PD(\< h \> \vsum U) \subset \PD(\KL \vsum U)=\PD(\EL). \]
\end{prop}
\begin{proof}
  By Example \ref{HolSymplPeriod} the period map is holomorphic.
  If $(h,e,f)$ is the standard basis of $\< h \> \vsum U$, then it is explicitly given by
  \[ p \mapsto [\exp(p h)] =[e \,+\, p h \,-\,\frac{1}{2}p^2 \Q{h}{h}
  f]  \in \PD(\< h \> \vsum U)\subset \PD(\EL).  \] 
  The injectivity of the period map is now obvious. To show surjectivity we
  let $[a e \,+\, b h \,+\, c f]$ be an arbitrary point in $\mathscr{D}(\< h \> \vsum U)$.
  By definition we have 
  \[     a c + 2 b^2 = 0, \quad   a \bar{c} + c \bar{a} + 4 b \bar{b} =  2 \Re(a \bar{c}) + 4 |b|^2  > 0 \]
  Hence $a \neq 0$ and we can set $p := b / a$. 
  Then $c/a =-2 p^2$, so that
  \[ [\exp(p h)] = [1 e \,+\, b/a h \,+\, c/a f] = [a e \,+\, b h
  \,+\, c f]. \]
  The inequality translates into $\Im(p)^2 > 0$. That means 
  \[ \IC \setminus \IR \lra \PD(\< h \> \vsum U),\; p \mapsto [\exp(p h)]   \] 
  is an isomorphism and therefore proves the proposition.
\end{proof}


\section{Period map for the Dwork family}\label{DworkSection}
\newcommand{\FF}{F}
\newcommand{\SF}{S}
\newcommand{\XF}{X}

\subsection{Construction of the Dwork family} \label{constr}
We start with the \emph{Fermat pencil} $\FF \subset \IP^3 \times \IP^1$ defined by the equation
\[ f = X_0^4 + X_1^4 + X_2^4 + X_3^4 - 4 t X_0 X_1 X_2 X_3 \]
where $X_0,\dots,X_3$ are homogeneous coordinates on $\IP^3$ and $t \in \IA^1 \subset \IP^1$ 
is an affine parameter.
We view $\FF$ as a family of quartics over $\IP^1$ via the projection $p:\FF \ra \IP^1$. 

The fibers $F_t=p^{-1}(\{t\})$ are smooth if $t$ does not lie in
\[ \Sigma = \{ t  \,|\, t^4=1 \} \union \{\infty \}. \]
For $t^4=1$ we find 16 singularities of type $A_1$,
for $t = \infty$ the Fermat pencil degenerates into the union of 
four planes: $X_0 X_1 X_2 X_3 = 0$.

Let $\mu_4$ denote the forth roots of unity. The group 
\[ G = \Set{ (a_0,a_1,a_2,a_3) }{  a_i \in \mu_4, a_0a_1a_2a_3 = 1 }/\mu_4 \isom (\IZ/4\IZ)^2 \]
acts on $\FF$ respecting the fibers $F_t$.

The quotient variety $\SF=\FF/G$ can be explicitly embedded into a 
projective space as follows. The monomials
\[ (Y_0, \dots, Y_4) :=(X_0^4,X_1^4,X_2^4,X_3^4,X_0 X_1 X_2 X_3)  \]
define a $G$-invariant map $\IP^3 \ra \IP^4$, and the image of $\FF$ in $\IP^4\times \IP^1$ under this 
morphism is cut out by the equations
\begin{equation}\label{SEQ}
  Y_0+Y_1+Y_2+Y_3-4 t Y_4,\quad Y_0 Y_1 Y_2 Y_3 - Y_4^4.
\end{equation}
It is easy to see that this image is isomorphic to the quotient $\SF$.

\begin{prop}
  For $t \neq \Sigma$ the space $S_t$ has precisely six singularities of type $A_3$.
  If $t^4=1$ there is an additional $A_1$-singularity.
  The fiber $S_\infty$ is a union of hyperplanes, it is in fact isomorphic to $F_\infty$ itself.
\end{prop}

\begin{proof}
  The first statement can be seen by direct calculation using (\ref{SEQ}).
  A more conceptual argument goes as follows.
  We note that the action of $G$ is free away from the $24$ points in
  \[  \union_{i\neq j} \Set{[X_0:\dots:X_3] }{ X_i = X_j = 0, F = 0 } \]
  which have stabilizer isomorphic to $\mu_4$. Around such a point $p$ we find an 
  analytic neighborhood $U$ such that the stabilizer $G_p$ acts on $U$ and $S_t$ 
  is locally isomorphic to $U/G_p$.
  \begin{center}\begin{minipage}{5cm} \xymatrix{
        F_t \ar[d]_{\phi} \ar@{<-^{)}}[r] &  U \ar[d] & \save \hspace{-1.8cm}\ni p \restore \\
        S_t \ar@{<-^{)}}[r]    &   U/G_p &
  }\end{minipage}\end{center}
  We can choose $U \subset \IC^2$ to be a ball on which $G_p \isom \mu_4$ acts
  as
  \[ a \cdot (x,y) = (a x, a^{-1} y) \]
  The quotient singularity is well known to be of type $D_2 =A_3$.

  To prove the second statement, recall that  there are $16$ singularities 
  of Type $A_1$ in each surface $F_t$ for $t^4=1$.
  It is easy to see that these form an orbit for the $G$ action and that they are disjoint from the
  $A_3$-singularities above.

  Finally, that $S_\infty$ is a union of hyperplanes follows directly form the equations (\ref{SEQ}).
\end{proof}

Note that $S_t \subset \IP^4$ is isomorphic to a (singular) quartic in $\IP^3$ since
the first equation defining $S_t$ is linear.

\begin{prop}
  There exists a minimal, simultaneous resolution of the $A_3$ singularities in $\SF \ra \IP^1$.
  That means, there is a threefold $\XF \ra \IP^1$ together with a morphism $\XF \ra \SF$
  over $\IP^1$ which restricts to a minimal resolution of the six $A_3$-singularities on each 
  fiber over $t \notin \Sigma$.
\end{prop}
\begin{proof}
  The position of the $A_3$-singularities of $S_t$ in $\IP^4$ does not change,
  as we vary $t$. So we can blow-up $\IP^4$ at these points. Also the singularities 
  of the strict transform of $S_t$ are independent of $t$. Hence we can construct
  $\XF$ by blowing-up the singularities again.
\end{proof}

\begin{defn} 
  The family $\XF \ra \IP^1$ is called the the \emph{Dwork Family}.

  The fibers $X_t$ are smooth for $t\in B = \IP^1 \setminus \Sigma$, 
  $\Sigma=\{ t \,|\, t^4 =1\} \union \infty$.  We denote by $\pi: \XF
  \ra B$ 
  the restriction.
\end{defn}

\begin{prop}
  The members $X_t$ of the Dwork family are $K3$ surfaces 
  for $t \notin \Sigma$.
\end{prop}
\begin{proof}
  It is shown in \cite{Nikulin1976}, that a minimal resolution of a
  quotient of a K3 surface by a finite group acting symplectically is
  again a K3 surface.
\end{proof}


\subsection{Holomorphic two-forms on the Dwork family}
In this subsection we construct holomorphic two-forms $\Omega_t$ on
the members of the Dwork family. We do this first for the Fermat
pencil using the residue construction (\cite{CMP2003} Section 3.3, 
\cite{GH1978} Chapter 5) and then pull back to the Dwork family.

Let $U := \IP^3 \setminus F_t $, there is a residue morphism:
\[ Res: \HH^k(U,\IC) \ra \HH^{k-1}(F_t,\IC). \]
This morphism is most easily described for de Rham cohomology groups.  
The boundary of a tubular neighborhood of $F_t$ in $\IP^3$ will be a 
$S^1$-bundle over $F_t$ completely contained in $U$. We integrate
a $k$-form on $U$ fiber-wise along this bundle to obtain a $k-1$ form on $F_t$,
this induces the residue map in cohomology.

\begin{rem}
  The residue morphism is also defined on the integral cohomology groups.
  It is the composition of the boundary morphism in the long exact sequence of 
  the space pair $(\IP^3,U)$ with the Thom isomorphism $\HH^{k+1}(\IP^3,U) \isom \HH^{k-1}(F_t)$ (up to a sign).
\end{rem}

There is a unique (up to scalar) holomorphic 3-form $\Xi_t$ on $\IP^3$, with simple poles along $F_t$. 
Its pull-back to $\IC^4\setminus \{0\}$ is given by the expression
\begin{align}\label{Omega}
  \Xi_t = \sum_{i=0}^3 (-1)^i \frac{X_i \, d X_0 \wedge \dots \wedge \widehat{dX_i} \wedge \dots \wedge d X_3}{X_0^4 + X_1^4 + X_2^4 + X_3^4 - 4 t X_0 X_1 X_2 X_3 }.
\end{align}
One checks that this form is closed and hence $\sigma_t:=Res(\Xi_t)$ is a well defined, closed two-form on $F_t$.

Let us choose coordinates $z_i=X_1/X_0, i=1,\dots,3$ for $\IP^3$, here
\[ \sigma_t=Res(\Xi_t)=Res(\frac{ d z_1 \wedge \dots \wedge  dz_3}{f_t } )\]
where $f_t =1 + z_1^4 + z_2^4 + z_3^4 - 4 t  z_1 z_2 z_3$ is the function defining $F_t$.

On the open subset $\del f_t / \del z_3 \neq 0$ the functions $(z_1,z_2)$ are (\'etale) coordinates for $F_t$, 
and $(f,z_1,z_2)$ are (\'etale) coordinates for $\IP^3$. In these coordinates the sphere bundle is just 
given by $|f|=\eps > 0$ and  fiber-wise integration reduces to taking the usual residue in each fiber $(z_1,z_2) = const$. 

Solving $df = \sum_i \del f / \del z_i dz_i $ for $dz_3$ and substituting above we get a local coordinate expression for $\sigma_t$:
\[ \sigma_t = Res(\frac{df}{f} \cdot \frac{dz^1 \wedge dz^2}{\del f / \del z_3}) = 2 \pi i  \frac{dz_1 \wedge dz_2}{ 4 z_3^3  - 4 t  z_1 z_2 }. \]

\begin{prop}
  The residue $\sigma_t$ of the meromorphic three form $\Xi_t$, 
  is a nowhere-vanishing holomorphic two form on all smooth 
  members $F_t$ of the Fermat pencil. \hfill $\qed$
\end{prop}

The same construction gives us a global version of $\sigma_t$: 
The inclusion $\FAM{F} \subset \IP^3 \times B$ is a smooth divisor, and the residue of the 
three-form $\Xi$ on $\IP^3 \times B$ given by same formula (\ref{Omega}) provides us with a two-form 
$\sigma$ on $\FAM{F}$ which defines a global section of $p_* \Omega^2_{\FAM{F}/B}$.
Clearly $\sigma$ restricts to $\sigma_t$ on each fiber and hence
trivializes the line bundle $p_* \Omega^2_{\FAM{F}/B}$.

We now proceed to the Dwork family. Consider the group 
\[ G = \Set{ (a_0,a_1,a_2,a_3) }{  a_i \in \mu_4, a_0 a_1 a_2 a_3 = 1 }/\mu_4 \]
acting on $F_t \subset \IP^3$. For  $g=(a_0,a_1,a_2,a_3) \in G$ we compute
\begin{align*}
g^*\Xi_t=\sum_{i=0}^3 (-1)^i \frac{a_0 \dots a_3 X_i \,X_0 \wedge \dots \wedge \widehat{dX_i} \wedge \dots \wedge d X_3}
{a_0^4 X_0^4 + a_1^4X_1^4 + a_2^4 X_2^4 + a_3^4 X_3^4 - 4 ta_0 \dots a_3 X_0 X_1 X_2 X_3 }
\end{align*}
which equals $\Xi_t$, hence $\sigma_t = Res(\Xi_t)$ is also $G$ invariant.
It follows that $\sigma_t$ descends to a form $\tilde{\sigma}$ on the smooth part 
$S_t^{reg} \subset S_t= F_t/G$.

Recall that the Dwork family is a simultaneous, minimal resolution of singularities 
$\rho: X_t \ra S_t$. In particular $\rho$ is an isomorphism over $S_t^{reg}$.

As $S_t^{reg}$ is isomorphic to an open subset of a K3 surface we find $\Omega^2_{S_t^{reg}} \isom \ko_{S_t^{reg}}$.
Moreover $\HH^0(S_t^{reg},\ko_{S_t^{reg}}) = \IC$ since the complement is an exceptional divisor.
It follows, that $\tilde{\sigma}$ extends to a holomorphic 2-form $\Omega_t$ on $X_t$.

The same construction works also in the global situation 
$\FF \ra \SF \leftarrow \XF$ over $B$
and gives us a global section $\Omega$ of $\pi_* \Omega^2_{\XF/B}$.

\begin{prop}
  There is a global section $\Omega$ of $\pi_* \Omega^2_{\FAM{X}/B}$ that restricts to $\Omega_t$ on each fiber.
  Moreover the pull-back of $\Omega_t$ along the rational map $\xymatrix{F_t \ar@{-->}[r] & X_t}$ coincides 
  with $\sigma_t$ on the set of definition.

  The section $\Omega$ trivializes the line bundle $\pi_* \Omega^2_{\FAM{X}/B}$
  and thus the variation of Hodge structures of $\pi: \XF \ra B$ is given by
  \[ \curly{F}^2 = \ko_B\; \Omega \;\subset\; \curly{F}^1 =
  (\curly{F}^2)^\perp \;\subset\; \HO \]
\end{prop}

\newcommand{\vspan}[1]{\negthinspace < \negmedspace #1 \negmedspace > \negthinspace }
\newcommand{\q}[1]{\<#1\>}
\newcommand{\blank}{\underline{\hspace{0.3cm}}}


\subsection{Monodromy of the Dwork family}
The Dwork family $\pi:\XF \ra B$ determines a local system
$\HOZ := \DR^2 \pi_* \IZ_\XF$ on $B$.
As is well known, every local system is completely determined by its
monodromy representation
\[ PT: \pi_1(B,t) \lra Aut((\HOZ)_t)=Aut(\HH^2(X_t,\IZ)), \quad t \in B \]
given by parallel transport. In this section we will explicitly
describe this representation.

To state the main result we need the following notation.  Let
\[ M_2=2E_8(-1) \oplus U \oplus  \<-4\>  \qtext{and} T_0 = \<4\> \oplus U. \]
If $e,f$ is the standard basis of $U$, and $l,h$ are generators of
$\<-4\>$ and $\<4\>$ respectively then we can define a primitive
embedding
\[ \<-4\> \oplus \<4\> \lra U,\quad l \mapsto e-2f,\; h \mapsto e+2f. \]
This induces also an  embedding $M_2 \oplus T_0 \ra 2 E_8(-1) \oplus 3 U = \KL$.

\begin{figure} \label{Paths}
  \caption{Generators of the fundamental group $\pi_1(B,t_0)$.}
  \includegraphics[width=60mm]{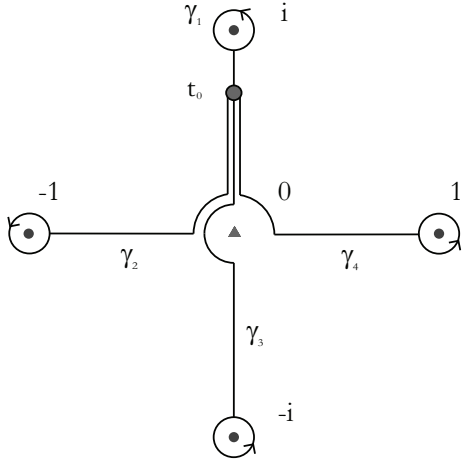}
\end{figure}

\begin{thm}[Narumiyah--Shiga, Dolgachev]\label{MonodromyThm}
  At the point $t_0=i/\sqrt{2}$ there is an isomorphism 
  \[ m: \HH^2(X_{t_0},\IZ) \lra \KL \]
  such that
  \begin{enumerate}
    \item[i)] The Neron--Severi group of each member $X_t$ contains the image of $M_2$ 
      under $m^{-1}$ composed with parallel transport along any path from $t_0$ to $t$ in $B$.
      For general $t$ this inclusion is an isomorphism.
    \item[ii)] The monodromy representation on $\HH^2(X_t,\IZ)$ respects the 
      images of the subspaces $M_2,T_0$ and acts trivially on the first one.
  \end{enumerate}
  Moreover, the monodromy representation on $T_0$ is given by the following matrices. 
  Let $(h,e,f)$ be the standard basis of $T_0=\<4\>\vsum U$, let $\gamma_k \in \pi_1(B,t_0)$ 
  be the paths depicted in Figure \ref{Paths} and  
  $\gamma_{\infty} = (\gamma_4\cdot \gamma_3\cdot \gamma_2\cdot
  \gamma_1)^{-1}$ the path around $\infty \in \IP^1$.\footnote{
  As Narumiyah and Shiga, we use the convention to compose paths like functions, 
  i.e. $\gamma: p \ra q,\, \delta: q \ra r$, then $\delta \cdot \gamma: p \ra r$. 
  This has the advantage, that monodromy becomes a representation, as opposed to 
  an anti-representation.}

  Then the following identities hold
  \[ PT_{\gamma_k} (h,e,f) = (h,e,f).M^k \]
  where
  \begin{gather*}  M^1=
  \begin{pmatrix} 
    1 & 0 & 0 \\
    0 & 0 & 1 \\ 
    0 & 1 & 0 
  \end{pmatrix},\quad 
   M^2=
  \begin{pmatrix} 
    5& 1 & -3 \\
  -12& -2 & 9 \\ 
    4&  1 & -2 
    \end{pmatrix}, \quad
    M^3=
  \begin{pmatrix} 
    17  & 6  &  -6 \\
    -24  & -8 &  9 \\ 
    24 & 9  & -8
  \end{pmatrix}\\
    M^4=
  \begin{pmatrix} 
    5   & 3 & -1 \\
    -4  & -2 & 1 \\ 
    12  & 9  &-2 
  \end{pmatrix}, \quad
    M^\infty=
  \begin{pmatrix} 
    1 & 4 &  0 \\
    0 & 1 &  0 \\ 
    -16 & -32  & 1
  \end{pmatrix}.
\end{gather*}
\end{thm}

\begin{rem}\label{mum_remark}
  Note that the matrix $M^{\infty}$ is unipotent of maximal index $3$, i.e.
  \[ (M^\infty-1)^3=0, \quad (M^\infty-1)^2 \neq 0 \]
  this will be crucial for the characterization of the period map in
  section \ref{sec:per_char}.
\end{rem}

The proof is a consequence of the following theorems.

\begin{thm}[Dolgachev \cite{Dolgachev1996}] \label{DolgachevTheorem}
  The Dwork family $\XF \ra B$ carries an $M_2$-polarization, 
  i.e.\ there exists a morphism of local systems
  \begin{equation} \label{M2P} Pol:  M_2 \tensor \underline{\IZ}_{B} \lra \HOZ   \end{equation}
  inducing a primitive lattice embedding in each fiber
  which factorizes through the inclusion $\Pic(X_t) \subset \HH^2(X_t,\IZ)$.
  Moreover, for general $t$ this map is an isomorphism onto the Neron--Severi group.
\end{thm}

\begin{thm}[Narumiyah--Shiga \cite{NS2001}] \label{NS_THM} 
  There is a primitive lattice embedding
  \[ Tr: T_0 \lra \HH^2(X_{t_0},\IZ) \]
  with image in the orthogonal complement of the polarization $Pol(M_2)_{t_0}$.
  Moreover the monodromy representation on $T_0$ is given by the matrices described in Theorem
  (\ref{MonodromyThm}).
\end{thm}

\begin{proof}
  The intersection form is stated in Theorem 4.1 of \cite{NS2001}, 
  and the monodromy matrices in Remark 4.2 following this theorem.
  We only explain how their notation differs form ours.

  They consider the family $\tilde{F}_\lambda \subset \IP^3$ defined by the equation
  \[X_0^4+X_1^4+X_2^4+X_3^4+\lambda  X_0 X_1 X_2 X_3 = 0.\]
  In order to ensure the relation $\lambda=4 t$ holds, we identify this family 
  via the isomorphism
  \[ \tilde{F}_\lambda \ra F_t, \quad X_0 \mapsto - X_0,\, X_1 \mapsto X_1,\, X_2 \mapsto X_2,\, X_3 \mapsto X_3 \]
  with our Fermat pencil.

  Their basis $(e',f',h')$ of $U \vsum \<4\>$ is related to our basis $(h,e,f)$ of $\< 4 \> \vsum U$ by
  \[ (h,e,f) = (e',f',h').T, \quad T:=\begin{pmatrix} 0 & 1 & 0 \\ 0 & 0 & 1 \\ 1 & 0 & 0 \end{pmatrix}  \]
  They introduce a new variable $t'=-(\lambda)^2/2$ and consider
  paths $\delta_1,\dots,\delta_3$ in the $t'$-plane (Fig.6 in \cite{NS2001}).
  The images of our paths $\gamma_1,\dots,\gamma_4$ are give by
  \begin{gather*} 
    \gamma_1 \mapsto \delta_1, \quad \gamma_2 \mapsto \delta_2^{-1} \cdot \delta_3 \cdot \delta_2, \quad
    \gamma_3 \mapsto \delta_2^{-1} \cdot \delta_1 \cdot \delta_2, \quad \gamma_4 \mapsto \delta_3.
  \end{gather*}
  Let $N_i$ be the monodromy matrices along $\delta_i$ as stated in Remark 4.2 of \cite{NS2001}.
  By what was said above, we compute the monodromy matrix e.g. along $\gamma_2$ as
  \[ M^2 = T^{-1}.(N_2)^{-1}.N_3.(N_2).T. \qedhere \]  
\end{proof}

So far we do not know whether the primitive embedding
\[ Pol_{t_0} \oplus Tr: M_2 \oplus T_0 \lra \HH^2(X_{t_0},\IZ)\] 
can be extended to an isomorphism of lattices $\KL \lra \HH^2(X_{t_0},\IZ)$.

The following theorem of Nikulin ensures, that we can always change $Pol_{t_0}$ by an 
automorphism of $M_2$ such, that an extension exists.

\begin{thm}[Nikulin \cite{Nikulin1979}, 1.14.4] \label{NIK_THM}
  Let $i:S \ra L$ be a primitive embedding of an even non-degenerate lattice $S$ of signature $(s_+,s_-)$
  into an even non-degenerate lattice of signature $(l_+,l_-)$.  For any other primitive embedding $j:S \ra L$,
  there is an automorphism $\alpha \in O(L)$ such that $i = j \circ \alpha$ if
  \[ l_+ > s_+,l_- > s_- \qtext{and} rk(L)-rk(S) \geq l(S) + 2 \]
  where $l(S)$ is the minimal number of generators of the discriminant group $S\dual/S$.
\end{thm}

We apply this theorem as follows. First choose an arbitrary 
isomorphism $\tilde{n}: \KL \ra \HH^2(X_{t_0},\IZ)$.
This gives us a primitive embedding of $T_0$ by restriction
\[ \tilde{n}|_{T_0} : T_0 \lra \HH^2(X_{t_0},\IZ). \]
Also there is the primitive embedding constructed in (\ref{NS_THM})
\[ Tr: T_0 \lra \HH^2(X_{t_0},\IZ). \]
Note that $sign(T_0)=(2,1),sign(\KL)=(3,19)$ and $l(T_0)=l(\<4\>)=1$,
so we can apply Nikulin's theorem to conclude, that these two
differ by an orthogonal automorphism $\alpha$ of $\HH^2(X_{t_0},\IZ)$.

Set $n= \alpha \circ \tilde{n}$ so that $\tilde{n}|_{T_0}=Tr$. 
Note also that $n$ induces an isomorphism of the orthogonal complements
\[ n|_{M_2}: M_2 = T_0^\perp \lra Tr(T_0)^\perp = Pol(M_2)_{t_0}. \]
As mentioned above, this isomorphism can differ, by an automorphism of $M_2$, from 
the one provided by Dolgachev's polarization.
It is now clear that $m = n^{-1}$ is a marking with the required properties. 
This concludes the proof of Theorem \ref{MonodromyThm}. \hfill $\qed$

\newcommand{\PO}{\mathcal{P}}
\newcommand{\POZ}{\PO_\IZ}
\newcommand{\POQ}{\PO_\IQ}
\newcommand{\POC}{\PO_\IC}

\newcommand{\TR}{\mathcal{T}}
\newcommand{\TRZ}{\TR_\IZ}
\newcommand{\TRQ}{\TR_\IQ}
\newcommand{\TRC}{\TR_\IC}

\begin{cor}\label{HOQ_DEC}
  The local system $\HO_\IQ := \HOZ \tensor \IQ$ 
  decomposes into an orthogonal direct sum
 \[  \HO_\IQ = \POQ \vsum \TRQ \]
  where $\POQ$ is a trivial local system of rank $19$ spanned by the 
  algebraic classes in the image of the polarization $Pol$,
  and $\TRQ$ is spanned by the image of $Tr$.
\end{cor}


\subsection{The Picard--Fuchs equation}
So far we have described the local system $\HOZ$ and the Hodge
filtration $\kf^i \subset \HO$ of the Dwork family independently.  The
next step is to relate them to each other by calculating the period
integrals
\[ t \mapsto \int_\Gamma \Omega_t \] 
for local sections $\Gamma \in \HOZ$.  The essential tool here is a
differential equation, the Picard--Fuchs equation, that is satisfied
by these period integrals.

Let $t$ be the affine coordinate on $B \subset \IA^1$, and $\del_t$
the associated global vector field.  The {\em Gau\ss--Manin
  connection} $\nabla$ on $\HO = \HOZ \tensor \ko_B$ is defined by
$\Gamma \tensor f \mapsto \Gamma \tensor df$.  We denote by
\[ \Omega^{(i)} := \nabla_{\del_t} \circ \dots \circ \nabla_{\del_t} \Omega \; \in \;  \HH^0(B,\HO) \]
be the $i$-th iterated Gau\ss--Manin derivative of $\Omega$ in direction $\del_t$.

\begin{prop}\label{PFRelation}
  The global section $\Omega$ of $\pi_*\Omega^2_{\XF/B}$ satisfies the 
  differential equation
  \begin{align} \label{PF} 
    \Omega^{(3)} = \frac{1}{1-t^4}(6 t^3 \Omega^{(2)} + 7 t^2 \Omega^{(1)} + t \Omega).
  \end{align}
\end{prop}

\begin{proof}
  This is an application of the Griffiths--Dwork reduction method,
  see \cite{Griffiths1969}, or \cite{Morrison1992} for a similar application.
  We will outline the basic steps.

  It is enough to prove the formula on the dense open subset $\rho(S_t^{reg})$ of $\XF$.
  Since the map $\FF \ra \SF$ is \'etale over $S_t^{reg}$, we can furthermore reduce
  the calculation to the Fermat pencil of quartic hypersurfaces.
  The holomorphic forms on the members $F_t$ of the Fermat pencil are 
  residues of meromorphic 3-forms $\Xi_t$  on $\IP^3$. 
  Since taking residues commutes with the Gau\ss--Manin connection, 
  we only need to differentiate the global 3-from $\Xi_t$. 

  We then use a criterion of Griffiths to show the corresponding equality
  between the residues. This involves a Gr\"obner basis computation 
  in the Jacobi ring of $F_t$. See e.g. \cite{smith2007} for an implementation.
\end{proof}

\newcommand{\DPF}{\curly{D}}

\begin{defn}
  We define the {\em Picard--Fuchs operator} associated to the Dwork
  family $\XF \ra B$ to be the differential operator
  \begin{align}\label{PicardFuchsEquation}
    \DPF = \del_t ^3 - \frac{1}{1-t^4}(6 t^3 \del_t^2 + 7 t^2 \del_t + t)
  \end{align}
  obtained from (\ref{PF}) by replacing $\nabla$ with
  $\del_t$.\footnote{See \cite{Morrison1992} or \cite{Peters1986} for
    a more general definition of the Picard--Fuchs equation.}
\end{defn}

\begin{rem}
  Let $\Gamma_t \in \HH^2(X_t,\IZ)$ be a cohomology class. 
  Extend $\Gamma_t$ to a flat local section $\Gamma$ of $\HOZ$.
  Since the quadratic form $\QF$ on $\HOZ$ is also flat, we can calculate
  \[  \del_t \<\Gamma,\Omega\> =\<\Gamma,\nabla_{\del_t} \Omega\> = \< \Gamma,  \Omega^{(1)} \>. \] 
  Similarly one finds that the function \[ t \mapsto \int_{\Gamma_t} \Omega_t = \<\Gamma,\Omega\>(t) \] 
  is a solution of the Picard--Fuchs equation $\DPF = 0$.
\end{rem}


\subsection{The period map of the Dwork family}
Recall that the Dwork family 
\[ \pi: \XF \lra B, \quad B = \IP^1 \setminus \Sigma, \quad \Sigma=\{ t \,|\, t^4 =1\} \union \infty\]
determines a variation of Hodge structures on $B$:
\[ \HOZ := R^2 \pi_* \underline{\IZ}_\XF, \quad 
\kf^2 = \pi_* \Omega^2_{\XF/B} \;\subset\;\kf^1 = (\kf^2)^\perp \;\subset\; \HO := \HOZ \tensor \ko_B. \]

We let $c: \tilde{B} \ra B$ be the universal cover, and choose a 
point $\tilde{t_0} \in \tilde{B}$ mapping to $t_0=i/\sqrt{2}$.

\begin{prop}
  The isomorphism constructed in Theorem \ref{MonodromyThm}
  \[ m: \HH^2(X_{t_0},\IZ) \lra \KL \]
  induces a marking of the local system $c^* \HOZ$. 
\end{prop}
\begin{proof}
  We compose $m$ with the canonical isomorphisms 
  \[ (c^*\HOZ)_{\tilde{t_0}} \lra (\HOZ)_{t_0} \lra \HH^2(X_{t_0},\IZ) \]
  and extend this map by parallel transport to an isomorphism of local systems
  \[ m:c^*\HOZ \lra \underline{\IZ}_{\tilde{B}} \tensor \KL.  \]
  This is possible since $\tilde{B}$ is simply connected, and hence
  both local systems are trivial.
\end{proof}

Choosing the marking in this way we get a period map
\[ \PER := \PER(c^*\curly{F}^*,m): \tilde{B} \lra  \PD(\KL). \]

\begin{prop} \label{OrthProp}
  Let $M_2,T_0 \subset \KL$ be as in Theorem \ref{MonodromyThm}.
  The period map takes values in $\PD(T_0) \subset \PD(\KL)$.
\end{prop}

\begin{proof}
  Let $D \in \HOZ$ be a local section contained in the 
  orthogonal complement $m^{-1}(M_2)$ of $m^{-1}(T_0)$.
  By Dolgachev's theorem \ref{DolgachevTheorem}, $D$ is fiber-wise contained in 
  the Picard group, hence $\Q{D}{\Omega} = 0$ by orthogonality of the Hodge decomposition.
\end{proof}

Let $(h,e,f)$ be the standard basis of $T_0 = \<4\>\oplus U$, we denote by the same symbols 
also the global sections of $c^*\HOZ$ associated via the marking.
By the last proposition we find holomorphic functions $a,b,c$ on $\tilde{B}$ such that
\begin{equation} c^* \Omega = a \, h + b \, e + c \, f  \;\in\, \HH^0(\tilde{B},c^* \HO) \end{equation}
and hence
\[ \PER =[a:b:c]: \tilde{B} \lra \IP(T_0 ) \subset \IP(\KL \tensor \IC),   \]
using the abusive notation  $[a:b:c]:=[a h + b e + cf]$.

\begin{rem}
  For each point $\tilde{p}\in \tilde{B},\; p=c(\tilde{p})$ there is a canonical isomorphism of stalks
  \[ c^*: \ko_{B,p} \lra \ko_{\tilde{B},\tilde{p}},\quad  f \mapsto f \circ c. \]
  In this way we may view functions on $\tilde{B}$ locally (on $\tilde{B}$) as functions on $B$.
\end{rem}

\begin{prop}
  If we view the functions $a,b,c$ locally as functions on $B$, then these functions satisfy 
  the Picard--Fuchs equation (\ref{PicardFuchsEquation}).
\end{prop}

\begin{proof}
  We can express $a,b,c$ as intersections with the dual basis in the following way. 
  If $(h\dual,e\dual,f\dual)=(h,e,f).G^{-1}=(1/4h,f,e)$, where $G$ is the Gram matrix 
  of $\QF$ on the basis $(h,e,f)$, i.e.
  \[  G=
  \begin{pmatrix} 
    4 & 0 & 0 \\
    0 & 0 & 1 \\ 
    0 & 1 & 0  
  \end{pmatrix}, \]
  then $a= \Q{h\dual}{\Omega}, b= \Q{e\dual}{\Omega}$ and $ c= \Q{f\dual}{\Omega}.$
  This exhibits the functions $a,b,c$ as period integrals and
  therefore shows that they satisfy the Picard--Fuchs equation.
\end{proof}

\begin{prop} \label{independence}
  The germs of the functions $a,b,c$ at $\tilde{p}$ form a basis for
  the three-dimensional vector space $Sol(\DPF,p) \subset \ko_{B,p}$
  of solutions of the Picard--Fuchs equation for all $\tilde{p} \in \tilde{B}$.
\end{prop}
\begin{proof}
  Linear independence of $a,b,c$ is equivalent to the non-vanishing of
  the Wronski determinant 
  \[ W = det \begin{pmatrix}
    a & b & c \\
    \del_t a & \del_t b & \del_t c  \\
    \del_t^2 a & \del_t b & \del_t^2  c 
  \end{pmatrix} \]
  of this sections. As the differential equation (\ref{PicardFuchsEquation})
  is normalized, this determinant is either identically zero or vanishes nowhere.\footnote{A standard 
    reference is \cite{Ince1944}, but see \cite{Beukers2007} for a readable summary.}
  If the vectors are everywhere linearly dependent, then we get a relation between
  the Gau\ss--Manin derivatives $\Omega, \Omega^{(1)}, \Omega^{(2)}$, since
  \[ \Omega^{(1)} = \nabla_{\del t} \Omega = (\del_t a) \, h + (\del_t b) \, e + (\del_t c) \, f. \]
  This means, that there is a order-two Picard--Fuchs equation for our family.
  That this is not the case, follows directly from the Griffiths--Dwork reduction 
  process (Proposition \ref{PFRelation}). 
\end{proof}

\subsection{Characterization of the period map via monodromies}\label{sec:per_char}
We have seen, that the coefficients of the period map satisfy the
Picard--Fuchs equation. In this section we characterize these
functions among all solutions. The key ingredient is the monodromy
calculation in Theorem \ref{MonodromyThm}.

\begin{rem}
  We briefly explain how analytic continuation on $B$ is related to
  global properties of the function on the universal cover $\tilde{B}$
  and thereby introduce some notation.

  Let $\tilde{p}$ be a point in $\tilde{B}$, mapping to
  $p=c(\tilde{p}) \in B$ and let $\delta: p \ra q$ be a path in
  $B$. There is a unique lift of $\delta$ to $\tilde{B}$ starting at
  $\tilde{p}$. Denote this path by $\tilde{\delta}: \tilde{p} \ra
  \tilde{q}$ and define $\delta \cdot \tilde{p} := \tilde{q}$.

  Also we can analytically continue holomorphic functions along
  $\delta$, this gives us a partially defined morphism between the
  stalks \[ AC_\delta: \ko_{B,p} \lra \ko_{B,q}.\] 
  A theorem of Cauchy \cite{Ince1944} ensures that if a
  function satisfies a differential equation of the form
  (\ref{PicardFuchsEquation}), then it can be analytically continued
  along every path.

  These two constructions are related as follows. Let $f:\tilde{B} \ra
  \IC$ be a holomorphic function. We can analytically continue the
  germ $f_{\tilde{p}} \in \ko_{B,p}$ along $\delta$ and get $
  AC_\delta f_{\tilde{p}} = f_{\tilde{q}}, \; \tilde{q}=\delta \cdot
  \tilde{p}$.
\end{rem}

Suppose now, that $\delta$ has the same start and end point
$t_0=i/\sqrt{2} \in B$.  We can express the analytic continuation of
$\PER$ along this paths in terms of the monodromy matrices of $\HOZ$.

\begin{prop}\label{PTMON}
  Let $\delta \in \pi_1(B,t_0)$ and 
  \[ PT_{\delta} (h,e,f) = (h,e,f).M^\delta   \]
  be the monodromy representation  of the local system $\HOZ$ as in Theorem \ref{MonodromyThm}.
  The analytic continuation of the period map at $\tilde{t_0}$ is given by
  \[ AC_{\delta} \PER_{\tilde{t_0}}= AC_{\delta} [a:b:c] = [a':b':c'] \]
  as tuple of germs at $\tilde{t_0}$, where
  \[ (a',b',c') = (a,b,c). G.M^\delta.G^{-1} \]
\end{prop}

\begin{proof}
  As remarked above we have the identity of tuples of functions on
  $\tilde{B}$
  \[ (a,b,c) = \Q{(h\dual,e\dual,f\dual)}{\Omega(p)}=\Q{\_}{\Omega} \circ (h,e,f).G^{-1}. \]
  Now integrals of the form $\int_{\Gamma(\tilde{p})}
  \Omega(p)=\Q{\Gamma}{\Omega}(\tilde{p})$ can be analytically
  continued by transporting the cycle $\Gamma$ in the local
  system. Thus we conclude
  \begin{align*} 
    AC_\delta (a,b,c) 
    &= AC_\delta \Q{\_}{\Omega} \circ (h,e,f).G^{-1} 
    = \Q{\_}{\Omega} \circ (PT_\delta(h,e,f)).G^{-1}\\
    &= \Q{\_}{\Omega} \circ (h,e,f).M^\delta.G^{-1}  
    = (a,b,c).G.M^\delta.G^{-1}. \qedhere
  \end{align*}
\end{proof}

\newcommand{\HG}[2]{ \hspace{0mm}_{#1}F_{#2} }
\newcommand{\FDD}{ \hspace{0mm}_3\curly{D}_2 }
\newcommand{\FD}{ \hspace{0mm}_2\curly{D}_1 }
\newcommand{\DT}{\DPF \circ t }
\newcommand{\PDP}{\tilde{exp}}
\newcommand{\IPDP}{\tilde{exp}^{-1}}
\newcommand{\PERC}{\PER^c}

We already saw in Proposition \ref{APeriodDomain} that the period
domain $\PD(\< 4 \> \vsum U)=\PD(T_0)$ is isomorphic to $\IC\setminus
\IR$. Let $(h,e,f)$ be the standard basis of $\< 4 \> \vsum
U$. A slightly different isomorphism is given by
\begin{align}\label{Param}
  \PDP:\IC\setminus \IR \lra \PD(T_0),\; z &\mapsto [z\,h\,-\,1\,e\,+\,2 z^2\,f] 
\end{align}
with inverse $\IPDP: [a\,h\,+\,b\,e\,+\,c\,f]  \mapsto -a/b$.

We consider the period map as a function to the complex numbers using
this parametrization of the period domain:
\[  \PERC = \IPDP \circ \PER: \tilde{B} \lra \IC. \]
We will see later, that the period map takes values in the upper half plane.

Theorem \ref{MonodromyThm} has a translation into properties of this function.
\begin{prop}
  The analytic continuation of the germ of the period map at $\tilde{t_0}$ along the 
  paths $\gamma_k$ depicted in Figure \ref{Paths} is given by
  \[ AC_{\gamma_k} \PERC_{\tilde{t_0}} = \beta_k(\PERC_{\tilde{t_0}})  \]
  where $\beta_k:\IH \ra \IH$ are the M\"obius transformations:
  \begin{gather*} 
    \beta_1(z)= \frac{-1}{2 z},\; \beta_2(z)= \frac{1 - 2 z}{2 - 6
      z},\; \beta_3(z)= \frac{3 - 4 z}{4 - 6 z}, \; \beta_4(z)= \frac{3 - 2 z}{2 - 2 z}, 
    \\  \beta_{\infty}(z)= 4 + z. 
  \end{gather*}
\end{prop} 
\begin{proof}
Direct calculation using Proposition \ref{PTMON}.
\end{proof}

The modification (\ref{Param}) of the parametrization was introduced to 
bring the monodromy at infinity to this standard form.

The fixed points of $\beta_i$ are 
\begin{gather}\label{FixedPoints}
  \beta_1:  \pm i/\sqrt{2}, \;\; 
  \beta_2: \frac{1}{3}(1 \pm i/ \sqrt{2}),  \;\;
  \beta_3: \frac{1}{3}(2 \pm i/ \sqrt{2}),  \;\;
  \beta_4: 1 \pm i / \sqrt{2}.
\end{gather}
These are also the limiting values of the period map at 
the corresponding boundary points $i,-1,-i,1 \in \IP^1\setminus B$.

The following characterization of the period map in terms of
monodromies is crucial. We show that the period map is determined up
to a constant by the monodromy at a maximal unipotent point
(cf. Remark \ref{mum_remark}). This is similar to the characterization
of the mirror map by Morrison \cite[Sec. 2]{Morrison1992}.  The remaining
constant can be fixed by considering an additional monodromy
transformation.

\begin{prop} \label{PER_CHAR}
  Let $a',b' \in \ko_{t_0}$ be non-zero solutions to the Picard--Fuchs equation 
  and $\PER':=a'/b'$. If
  \begin{gather*} 
    AC_{\gamma_\infty} (a',b') = (a',b').\begin{pmatrix} 1 & 0 \\ 4 & 1 \end{pmatrix}, \quad 
  \end{gather*}
  then there is a $\mu \in \IC$ such that $\PER' = \PERC + \mu$ as germs at $\tilde{t_0}$.

  If furthermore \[ AC_{\gamma_1} \PER' = \beta(\PER')\]  for a M\"obius 
  transformation $\beta$ with fixed points $\pm i/\sqrt{2}$, then $\PER' = \PERC$.
\end{prop}

\begin{proof}
  By Proposition \ref{independence} the functions $a',b'$ are a $\IC$-linear 
  combination of $a,b,c$.
  The monodromy transformation of $(a,b,c)$ at infinity is
  \[ N^{\infty}:=G.M^{\infty}.G^{-1}=
    \begin{pmatrix}
      1& 0& 16 \\
      -4& 1& -32 \\
      0& 0& 1
  \end{pmatrix}. \]
  Note that $a,b$ have the same monodromy behavior as $-a',b'$ at infinity.
  The matrix $N^\infty$ is unipotent of index $3$, i.e. $(N^{\infty}-\id)^3 = 0, \,(N^{\infty}-\id)^2 \neq 0$.
  In particular the only eigenvalue is $1$ and the corresponding eigenspace is one-dimensional,
  spanned by $e_2=(0,1,0)^t$. Hence there is a $\lambda \in \IC$ such that $b' = \lambda b$.

  The vector $v=(1,0,0)^t$ is characterized by the property $(N^{\infty}-1)v=-4e_2$.
  The space of such $v$ is a one dimensional affine space over the eigenspace $\IC\, e_2$.
  We conclude that $-a' =  \lambda a - \mu b$, for some $\mu \in \IC$.
  Since $b \neq 0$ it is $\lambda \neq 0$ and we may assume $\lambda = 1$.
  Hence 
  \[ \PER'=a'/b' = -a/b + \mu = \PERC + \mu. \]
  Moreover the monodromy of this function along $\gamma_1$ is
  \[ AC_{\gamma_1} \PER' = AC_{\gamma_1} \PERC + \mu = \beta_1(\PERC) + \mu. \]
  The fixed point equation $\beta_1(z) + \mu=z$ is a polynomial of degree $2$
  with discriminant $-2 + \mu^2$. This means the difference of the two 
  solution is $i \sqrt{2}$ only if $\mu=0$.
\end{proof}

\subsection{Nagura and Sugiyama's solutions}\label{NSSols}
Solutions to the Picard--Fuchs equation matching the criterion
\ref{PER_CHAR} were produced by Nagura and Sugiyama in
\cite{NS1995}. To state their result, we first need to transform the
equation.

The first step is to change the form $\Omega$ to $t^{-1} \Omega$,
which does not affect the period map, but changes the Picard--Fuchs
equation from $\DPF=0$ to $\DPF.t=0$.  We can further multiply by
$(1-t^4)$ from the left, without changing the solution space.  This
differential equation now does descend along the covering map
\[ z:B\setminus \{ 0 \} \lra \IP^1 \setminus \{0,1,\infty\}, \quad t \mapsto z(t)=t^{-4} \]
to a hypergeometric system on $C$. 
\begin{prop}\label{PFtoHG}
  Let 
  \begin{align}\label{F32}
    \FDD := \vartheta^3-z (\vartheta+1/4) (\vartheta +2/4) (\vartheta +3/4), \quad \vartheta = z \del_z
  \end{align}
  be the differential operator on $\IP^1 \setminus \{0,1,\infty\}$ associated to the generalized 
  hypergeometric function $\HG{3}{2}(1/4,2/4,3/4;1;1;u)$ then
  \[ z^* \FDD = \frac{1}{64}(1-t^4).\DPF.t. \]
\end{prop}
\begin{proof} Direct calculation.
\end{proof}

\begin{ex}\label{ExHGF}
  The function on $B$
  \[ \HG{3}{2}(\frac{1}{4},\frac{2}{4},\frac{3}{4};1;1;t^{-4}) t \]
  defined for $|t|>1$ satisfies the Picard--Fuchs equation.
\end{ex}

Consider the solutions to the hypergeometric differential equation $\FDD$
\begin{align*}
  W_1(z) &= \sum_{n=0}^{\infty} \frac{(4n)!}{(n!)^4 (4^4)^n}z^{n} = \HG{3}{2}(\frac{1}{4},\frac{2}{4},\frac{3}{4};1,1;z) \\
  W_2(z) &= \ln(4^{-4} z) W_1(z)  + 4 \sum_{n=0}^{\infty} \frac{(4n)!}{(n!)^4 (4^4)^n}[\Psi(4n+1)-\Psi(n+1)]z^{n}.
\end{align*}
where $\Psi$ denotes the digamma-function $\Psi(z)=\Gamma'(z)/\Gamma(z)$. 
The functions $W_i(t^{-4}), i=1,2$ are solutions to the pulled back
equation $z^*\FDD$. We set 
\[   P(t) := \frac{1}{2\pi i} \frac{W_2(t^{-4})}{W_1(t^{-4})}. \]

These functions converge for $|t|>1$ and hence define germs at the point
$t_1 =i\sqrt{2}$. The logarithm is chosen in such a way that $\Im(\ln((4 t_1)^{-4}))=0$.

Choose a path $\delta: t_0 \ra t_1,\, t_0=i/\sqrt{2}, t_1=i \sqrt{2}$ within  
the contractible region $\{t\,|\,\Re(t)>0,\, \Im(t)>0\} \subset B$.
We get an isomorphism between the fundamental groups by
\[  T_\delta: \pi_1(B,t_0) \ra \pi_1(B_,t_1), \quad \gamma \mapsto  \delta \cdot  \gamma \cdot \delta^{-1} \]
The analytic continuation along $T_\delta \gamma_{\infty}$ can be read off the definition
\[ AC_{T_\delta \gamma_{\infty}} W_1 = W_1, \quad AC_{T_\delta \gamma_{\infty}} W_2 = W_2 + 4 (2 \pi i) W_1. \]
Indeed, the sums define holomorphic functions and are therefore
unaffected by analytic continuation. The only contribution comes 
from the logarithmic term.
The path $T_\delta \gamma_{\infty}$ encircles $\infty$ once with positive orientation.
Therefore $0$ is encircled with negative orientation, so the logarithm picks up a summand $-2 \pi i$.

We can apply the first part of criterion  \ref{PER_CHAR} to see
\[ \PERC(t) = P(t) + \mu \]
as germs of functions at $\tilde{t_1}:=\delta \cdot \tilde{t_0}$ for some $\mu \in \IC$. 
To apply the second part of the criterion we need the following
additional information.

\begin{thm}[Nagura, Sugiyama \cite{NS1995}] \label{GlobalMon}
  An analytic continuation of the map $P:=\frac{1}{2 \pi i} W_2 / W_1$ to a sliced neighborhood of $t=1$
  is given by
  \begin{align*} 
    P(t)   &= \frac{i}{\sqrt{2}} \frac{U_1(t) + \tan(\frac{\pi}{8})^{-1} U_2(t) }{U_1(t) - \tan(\frac{\pi}{8})^{-1} U_2(t)} \\
    U_1(t) &= \frac{\Gamma(\frac{1}{8})^2}{\Gamma(\frac{1}{2})} \HG{2}{1}(\frac{1}{8},\frac{3}{8};\frac{1}{2};1-t^4) \\
    U_2(t) &= \frac{\Gamma(\frac{5}{8})^2}{\Gamma(\frac{3}{2})} (t^4 - 1)^{1/2} \HG{2}{1}(\frac{5}{8},\frac{5}{8};\frac{3}{2};1-t^4).
  \end{align*}
  Thus the monodromy around the point $t=1$ satisfies $ AC_{T_\delta \gamma_4} P  = -\frac{1}{2 P}$.
\end{thm}

We find get following corollary.

\begin{thm} \label{mainthm}
  The composition of the period map with the parametrization of the period domain (\ref{Param})
  \[ \PERC = \PDP \circ \PER: \tilde{B} \lra \PD(T_0) \lra \IC \setminus \IR \] 
  is explicitly given in a neighborhood of $\tilde{t_1}$ by 
  \[ \PERC(t) = P(t) = \frac{1}{2 \pi i} \frac{W_2(t^{-4})}{W_1(t^{-4})}. \]
\end{thm}

\begin{proof}
   We have to check, that the function $P$ has the right analytic 
   continuation along $T_\delta \gamma_1$, i.e. $ AC_{ T_{\delta} \gamma_1} P = -1/(2 P)$. 
   We know the analytic continuation of $P$ along $T_{\delta} \gamma_4$ has this form.

   But $P$ only depends on $z=t^{-4}$ not on $t$ itself. Moreover
   the images of the paths $\gamma_1$ and $\gamma_4$ under $t \mapsto t^{-4}$
   coincide.  Hence also the analytic continuations are the same.
\end{proof}

\begin{prop}\label{MMExpansion}
  The power series expansion of $\PERC$ at $t=\infty$ is given by
  \[ \PERC(w) = \frac{1}{2 \pi i } ( \ln(w)+104\,w+9780\,{w}^{2}+{4141760}/{3}\,{w}^{3}+231052570\,{w}^{4}+ \dots) \]
  \[ exp(2\pi i \PERC(w)) = w+104\,{w}^{2}+15188\,{w}^{3}+2585184\,{w}^{4}+480222434\,{w}^{5} +\dots \]
  where $w=1/(4 t)^{4}=4^{-4} z$.
\end{prop}

This is precisely the series obtained by Lian and Yau \cite{LY1996}\footnote{
  Equation 5.18 contains an expansion of the inverse series to ours.}
using a different method (see Remark \ref{Comparison}). They also prove that the 
expansion of $exp(2\pi i \PERC(w))$ has integral coefficients.

\begin{cor}
  The period map $\PERC$ takes values in the upper half plane.
\end{cor}

 \begin{rem} \label{PeriodThmProof}
   We show how Theorem \ref{PeriodThm} stated in the introduction can be derived from \ref{mainthm}.
  
  We identify $\HH^2(X_{t_0},\IZ) \isom \Lambda$ via the isomorphism given in Theorem \ref{MonodromyThm}
  and use parallel transport to extend this isomorphism to nearby fibers $X_t$.

  The period vector $\Omega_t$ is contained in $T_0 \tensor \IC$, 
  where $T_0=\<4\> \oplus U \subset \Lambda$ is the generic transcendental lattice. 
  By Theorem \ref{mainthm} and (\ref{Param}) 
  we have \[ [ \Omega_t ] = [ \PDP(\PERC(t)) ] \in \kd(T_0) \subset \IP(\Lambda_\IC) \] 
  and hence there is a nowhere vanishing holomorphic  function $f(t)$ such that 
  \begin{align} \label{PerCor}  f(t) \Omega_t=\PDP(\PERC(t))=\PERC(t) h - e +2 (\PERC(t))^2 f. 
  \end{align}
  As $f(t) \Omega_t$ is also a non-vanishing holomorphic two-form we
  can assume this equation holds true already for $\Omega_t$.  The
  period integrals can now be calculated as intersection products
  $\int_\Gamma \Omega_t = \Omega_t . \Gamma$.
  
  The required basis $\Gamma_i$ of $\Lambda = 2E_8(-1) \oplus U'' \oplus
  U' \oplus U$ is constructed as follows.  We let
  $(\Gamma_1,\Gamma_2,\Gamma_3)=(h,e,f)$ be the standard basis of
  $T_0$.  Recall that $h=e' + 2f'$ and hence
  $(\Gamma_1,\Gamma_4)=(h,f')$ is a basis of $U'$.  The remaining
  basis vectors can be chosen to be any basis of the orthogonal
  complement $2E_8(-1) \oplus U''$ of $(\Gamma_1,\dots,\Gamma_4)$.  Using
  (\ref{PerCor}) it is now straightforward to calculate the entries of
  the period vector.
 \end{rem}


\subsection{The period map as Schwarz triangle function}

In this section we will relate the period map to a Schwarz triangle function.
We begin by recalling some basic facts about these functions from \cite{Beukers2007}.

\begin{defn}
  The {\em hypergeometric differential equation} with parameters $a,b,c \in \IC$ is
  \begin{align}
    \vartheta(\vartheta +c -1) f -z(\vartheta + a) (\vartheta + b) f = 0, \quad \vartheta = z\del_z,\, f\in \ko_\IC
  \end{align}
  which is satisfied by the hypergeometric function $f=\HG{2}{1}(a,b;c;z)$.
  
  Let $f,g$ be two independent solutions to this differential equation at a point $z_0 \in \IH$.
  The function $D(z)=f/g$ considered as map $\IH \ra \IC$ is called \em{Schwarz triangle function}.
\end{defn}

These functions have very remarkable properties and were studied extensively 
in the 19th century (see Klein's lectures \cite{KleinHYP}).

\begin{defn}
  A {\em curvilinear triangle} is an open subset of $\IP^1$ whose boundary is the union of
  three open segments of  circles or lines and three points.
  The segments are called {\em edges} and the points {\em vertices} of the triangle.
\end{defn}

\begin{prop}\label{ExTr}
  For any three distinct points $A,B,C \in \IP^1$ and positive, real numbers $\lambda,\mu,\nu$ with
  $\lambda + \mu + \nu < 1$ there is a unique curvilinear triangle with
  vertices $(A,B,C)$ and interior angles $(\lambda \pi,\mu \pi, \nu \pi)$ in that order.
\end{prop}

\begin{thm}[Schwarz, \cite{Beukers2007} 3.20]
  A Schwarz triangle function maps the closed upper half plane $\IH \union \IR$ isomorphically
  to a curvilinear triangle. 
  
  The vertices are the points $(D(0),D(1),D(\infty))$ and the corresponding angles $(\lambda \pi,\mu \pi, \nu \pi)$ 
  depend on the parameters of the hypergeometric differential equation via $\lambda=|1-c|,\mu=|c-a-b|,\nu=|a-b|$.
\end{thm}

\newcommand{\TILDE}{\hspace{-1mm}\widetilde{\phantom{X}}}
\newcommand{\UC}{(\IP^1\setminus\{0,1,\infty\})\hspace{-1mm}\widetilde{\phantom{X}}}
\newcommand{\IPP}{\IP^1\setminus\{0,1,\infty\}}

Recall that the period map is a function on the universal cover 
of $B=\IP^1 \setminus \Sigma$ to the upper half plane.
\[ \PER^c: \tilde{B} \lra \IH. \]
This maps  descends along $t \mapsto z(t)=t^{-4}$ to a multi-valued map 
on $\IP^1 \setminus \{0,1,\infty\}$. We explain this last sentence more formally.
The map $t \mapsto z(t)=t^{-4}$ is an unramified covering $B \setminus \{0\} \ra \IP^1\setminus\{0,1,\infty\}$.
Hence it induces an isomorphism between the universal covering spaces.
Moreover the inclusion $B \setminus \{0\}\ra B$ induces a map 
$(B\setminus \{0\})\TILDE \ra \tilde{B}$. 
We use the composition
\[  \UC \isom (B\setminus \{0\})\TILDE \lra \tilde{B} \]
to view $\PER^c:\tilde{B} \ra \IH$ as multi-valued map on $\IP^1\setminus\{0,1,\infty\}$.

We choose a basepoint $\tilde{z}_1$ of $\UC$ mapping to $\tilde{t_1}$. Denote by  $\iota$ the unique lift of the 
inclusion $\IH \ra \IPP$ to the universal cover of $\IPP$ mapping $z_1$ to $\tilde{z}_1$ (when extended to the boundary of 
$\IH \subset B$).

\begin{thm}\label{TriangleThm}
  The restriction of the period map 
  \[ \PER^c(z):\UC \lra \IP^1 \]
  to $\iota:\IH \ra \UC$ is a Schwarz triangle function.
  The upper half plane is mapped to the triangle with vertices $(\infty,\frac{i}{2},\frac{1+i}{2})$
  and angles $(0,\pi/2,\pi/4)$ as pictured in Figure \ref{FigureIntro} in the introduction.
\end{thm}

\begin{proof}
  The strategy is the following. 
  We first construct the a triangle function with the expected mapping behavior.
  Then we write this function as a quotient of solution of the Picard--Fuchs equation.
  Finally we show that the assumptions of Proposition \ref{PER_CHAR} are satisfied by this 
  function. It follows that it has to be the period map.

  {\em Step 1.}
  Let $f,g$ be two independent solutions to $\FD$ at $t_1$. 
  By Schwarz' theorem $D(z)=f/g$ is a triangle function. 
  Using a M\"obius transformation, we can change the vertices of the triangle to be $(0,1,\infty)$. 
  As the composition is again of the form $f'/g'$ for independent solutions $f',g'$ of $\FD$
  we can assume $D(z)$ maps $(0,1,\infty)$ to $(\infty,\frac{i}{\sqrt{2}},\frac{1+i}{2})$.

  The triangle pictured in green color in Figure \ref{FigureIntro} 
  is the unique curvilinear triangle with vertices $(\infty,\frac{i}{\sqrt{2}},\frac{1+i}{2})$ 
  and interior angles $(0,\pi/2,\pi/4)$. Hence it is the image of $\IH$ under $D(z)$.

  The analytic continuation of $D(z)$ can be obtained by reflecting
  the triangle at its edges. This technique is called 
  {\em Schwarz reflection principle} (see \cite{Beukers2007} for details).

  Let $\delta_0,\delta_1 \in \pi_1(\IP^1\setminus\{0,1,\infty\},z_1)$
  be the paths pictured in Figure \ref{LeveltPaths} encircling $0,1$ once 
  with positive orientation respectively.
  Reflecting the triangles according to the crossings of the paths 
  with the components of $\IR\setminus\{0,1\}$ we find
  \[ AC_{\delta_0} D(z) =  D(z) + 1, \quad AC_{\delta_1} D(z) =  \frac{-1}{2D(z)}. \]
  
  This means that $AC_{\delta_0} (f/g) = (f+g)/g$ and 
  since $f,g$ are independent we can conclude that there is a $\lambda \in \IC^*$ such that
  \begin{align}  \label{MonMatrix}
    AC_{\delta_0} (f,g)= (f,g).
    \begin{pmatrix} \lambda & 0 \\ \lambda & \lambda
    \end{pmatrix}. 
  \end{align}
  The hypergeometric function $\HG{2}{1}(\frac{1}{8},\frac{3}{8};1;z)$ is a linear 
  combination of the basis solutions $(f,g)$. 
  Since it is holomorphic at $0$, the matrix (\ref{MonMatrix}) has to 
  have the eigenvalue $1$ which is only the case if $\lambda = 1$.

  \begin{figure}\label{LeveltPaths}
    \includegraphics[width=5cm]{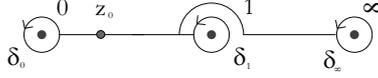}
    \caption{The paths $\delta_i$ in $\IP^1\setminus\{ 0,1,\infty \}$ based at $z_1=1/4$. }
  \end{figure}

  {\em Step 2.}
  The $\HG{3}{2}$-hypergeometric function $W_1(z)$ occuring in the expansion of 
  the period map is related to a $\HG{2}{1}$-hypergeometric function 
  by the Clausen identity (\cite{Bailey}, p.86)
  \[  \HG{3}{2}(\frac{1}{4},\frac{2}{4},\frac{3}{4};1;1;z) =\HG{2}{1}(\frac{1}{8},\frac{3}{8};z)^2.  \]
  The corresponding statement in terms of differential equations reads as follows.

  \begin{prop}
    The differential equation
    \begin{align}\label{GaussHG}
      \FD = \vartheta^2-z(\vartheta + 1/8) (\vartheta + 3/8), \quad \vartheta = z \frac{\del}{\del z}
    \end{align}
    associated to the hypergeometric function $\HG{1}{2}(\frac{1}{8},\frac{3}{8};1;z)$
    has the property that for all solutions $f,g$ to $\FD$ the product satisfies 
    $\FDD ( f.g ) = 0$.
    
    Conversely any solution to $\FDD$ is a sum of products of solutions to $\FD$.
  \end{prop}
  \begin{proof}
    The proposition can be rephrased by saying $\FDD = Sym^2(\FD)$.
    There is an algorithm to compute such symmetric squares of differential operators, 
    which is implemented e.g. in Maple. We used this program to verify the equality.
  \end{proof}

  Using this proposition and Proposition \ref{PFtoHG} 
  we can trivially express $D(z)$ as a quotient 
  of solutions of the Picard--Fuchs 
  equation (\ref{PicardFuchsEquation}),  namely 
  \[ D(t^{-4}) = \frac{f(t^{-4})}{g(t^{-4})} = \frac{f(t^{-4}) g(t^{-4})  \, t}{g(t^{-4})^2\, t}. \]

  {\em Step 3.}
  We claim that the tuple $(a,b)=(f(t^{-4}) g(t^{-4})  \, t,g(t^{-4})^2\, t)$
  of solutions of the Picard--Fuchs equation satisfies the assumptions of the 
  criterion \ref{PER_CHAR}.
  
  The paths $T_\delta \gamma_\infty, T_\delta \gamma_1 \in \pi_1(B,t_1)$ in $B$ 
  map to $\delta_0^4,\delta_1 \in \pi_1(\IPP,z_1)$ under $t \mapsto z(t)=t^{-4}$.  
  Hence we can calculate the monodromy transformations as
  \[ 
  AC_{T_\delta \gamma_\infty} (f,g) = (f,g). \begin{pmatrix} 1 & 0 \\ 1 & 1 \end{pmatrix}^4 =(f,g).\begin{pmatrix} 1 & 0 \\ 4 & 1 \end{pmatrix} \]
  and consequently also
  \[ 
  AC_{T_\delta \gamma_\infty} (a,b)  = (a,b). \begin{pmatrix} 1 & 0 \\ 4 & 1 \end{pmatrix} 
  \]
  moreover
  \[ AC_{T_\delta \gamma_1}  D(t^{-4}) = \frac{-1}{2D(t^{-4})} \]
  as required. This concludes the proof of the theorem.
\end{proof}


\section{Mirror symmetries and mirror maps}\label{MirroMapChapter}
It remains to translate the above computations in the framework
developed in chapter \ref{Framework}.

Let $\FAM{X} \ra B$ be the Dwork Pencil and
\[ \PER_B: \tilde{B} \lra \PD(T_0) \subset \PD(\KL) \subset
\PD(\EL) \] the (B-model) period map associated to the marking,
constructed in Theorem \ref{MonodromyThm}.  Here $T_0 \isom \< h \>
\vsum U$ is the transcendental lattice of the general member of
$\FAM{X}/B$.

Let $\curly{Y} \ra \IH$ be the family of generalized K3 structures on 
a quartic $Y \subset \IP^3$ as constructed in section \ref{QuarticSection} and
\[ \PER_A: \IH \lra \PD(\< H \> \vsum U) \subset \PD(\EL)  \]
the A-model period map as in Proposition \ref{APeriodDomain}. Here $\< H \> \vsum U$ is the lattice 
spanned by the class of a hyperplane $H$ and $U \isom \HH^0 \vsum \HH^4 \subset \HT(Y,\IZ)$.

\begin{thm}\label{MMTheorem}
  Mirror symmetry as described in Section \ref{MSFamily}
  between the symplectic quartic in $\IP^3$ and the Dwork family  is determined 
  by the diagram
  \begin{center}
    \begin{minipage}{5cm} 
      \xymatrix{
        \tilde{B} \ar[rr]^{\PER_B} \ar[d]_{\psi} & & \PD(T_0) \ar[r] \ar[d]^{g_0} & \PD(\EL)\ar[d]^g   \\
        \IH \ar[rr]^{\PER_A } & & \PD(\< H \> \vsum U) \ar[r]  &  \PD(\EL) 
      }
    \end{minipage}
  \end{center}
  where $g \in \ko(\EL)$ is a isometry interchanging $\HH^0 \vsum \HH^4$ with $U \subset T_0$
  and $\psi = \PER^c$ is the period map of Theorem \ref{mainthm}.
\end{thm}

\begin{proof}
  Recall from \ref{APeriodDomain} that $\PER_A(z)=[1e + z H -2 z^2 f]$.
  On the other hand $\PER^c$ was defined using the parametrization
  $\PDP(z)=[-1 e + zh + 2z^2 f]$. So in order for the diagram to commute
  we should use the isometry
  \[ g_0 :  T_0 = (\< h \> \vsum U) \lra (\< H \> \vsum U), \; h \mapsto H, e \mapsto -e, f \mapsto -f \]
  to relate the period domains $\PD(T_0)$ and $\PD((\< H \> \vsum U))$.
  This isomorphism is easily seen to extend to an isometry g of $\EL$ using 
  Nikulin's theorem \ref{NIK_THM}.
\end{proof}

\begin{rem}\label{Comparison}
  A period map in the sense of Morrison \cite{Morrison1992} is a 
  quotient $\psi=a/b$ of two solutions to the Picard--Fuchs 
  equation $a,b$ satisfying the property
  \[ AC_{\gamma_\infty} \psi = \psi + 1 \] for analytic continuation
  around the point of maximal unipotent monodromy.  As in Proposition
  \ref{PER_CHAR} one finds that $\psi$ is uniquely determined up to
  addition of a constant.  One chooses this constant in such a way
  that the Fourier expansion at $\infty$ has integral coefficients.

  Such a function can be constructed directly from the differential equation
  by using a {\em Frobenius basis} for the solutions at the singular point.
  Using this method, Lian and Yau \cite{LY1996} arrive at precisely the same 
  formula \ref{MMExpansion}.  

  There are several differences to our definition. 
  First note, that our mirror maps are symmetries of the period domain
  of (generalized) K3 surfaces which become functions only after composition
  with the corresponding period maps.

  Secondly and more importantly, we do require the solutions $a,b$ to be
  of the form $\int_\Gamma \omega_t$, for some {\em integral} cycle $\Gamma \in 
  \HH^2(X_{t_0},\IZ)$. It is not clear (and in general not true) that 
  the Frobenius basis has this property. 
  This was the main difficulty we faced above. Our solution relied heavily 
  on the work of Narumiyah and Shiga \cite{NS2001}.

  There is also a conceptual explanation that Morrison's mirror map
  coincides with ours. Conjecturally (see \cite{KKP2008}, \cite{Iritani2009}) the Frobenius 
  solutions differ form the integral periods by multiplication with the $\hat{\Gamma}$-class 
  \[ \hat{\Gamma}(X) = \prod_{i=1}^n \Gamma(1+\delta_i) = exp(-\gamma c_1(X) + \sum_{k \geq 2} (-1)^k (k-1)! \zeta(k) ch_k(TX)) \]
  where $\delta_i$ are the Chern roots of $TX$, $\gamma$ is Euler's constant and $\zeta(s)$ is the Riemann zeta function.
  The Calabi--Yau condition $c_1(X)=0$ translates into the statement, that the first two 
  entries of the Frobenius basis give indeed integral periods. 
  In our case, this information suffices to fix the Hodge structure completely.
\end{rem}


\bibliographystyle{alpha}
\bibliography{myReferences}

\end{document}